\documentclass[12pt]{article}

\usepackage{amsmath,amsfonts,amssymb,amsthm}

 \usepackage{ragged2e}

\def \Z{\mathbb{Z}}

\def \N{\mathbb{N}}
\def \C{\mathbb{C}}

\def \k{\bf k}
\def \i{\bf i}
\def \j{\bf j}
\def \a{\bf a}
\def \b{\bf b}
\def \c{\bf c}

\def \p{\prime}

\def \es{\epsilon}

\def \M{\mathbb{M}}
\newcommand{\CL}{{\mathcal L}}

\def \L{{\mathcal L}_{\alpha,\beta,F}}
\def \tl{\tilde{\mathcal{L}}}
\def \ud{\underline{d}}

\newcommand{\vf}{\varphi}
 \def \ad{\text{\rm ad}}

\def \e {\ell_{123}}

\def \id {\mbox{Ind}}
\def \d {\mbox{deg}}

\def \be{\begin{equation}\label}
\def \ee{\end{equation}}
\def \bex{\begin{example}\label}
\def \eex{\end{example}}
\def \bl{\begin{lem}\label}
\def \el{\end{lem}}
\def \bt{\begin{thm}\label}
\def \et{\end{thm}}
\def \bp{\begin{prop}\label}
\def \ep{\end{prop}}
\def \br{\begin{rem}\label}
\def \er{\end{rem}}
\def \bc{\begin{coro}\label}
\def \ec{\end{coro}}
\def \bd{\begin{de}\label}
\def \ed{\end{de}}

\newtheorem{thm}{Theorem}[section]
\newtheorem{prop}[thm]{Proposition}
\newtheorem{coro}[thm]{Corollary}

\newtheorem{example}[thm]{Example}
\newtheorem{lem}[thm]{Lemma}
\newtheorem{rem}[thm]{Remark}
\newtheorem{de}[thm]{Definition}

\input amssym.def
\begin{document}

\begin{center}
{\Large \bf  Restricted modules and associated vertex algebras of extended Heisenberg-Virasoro algebra}
\end{center}

\begin{center}
	{Hongyan Guo
		\footnote{Partially supported by China NSF grants (No.11901224)
and the Fundamental Research Funds for
the Central Universities (CCNU22QN002)},
 Huaimin Li
		\\
School of Mathematics and Statistics,
and Hubei Key Laboratory of Mathematical Sciences,
Central China Normal University, Wuhan 430079, P. R. China
	}
	
\end{center}

\begin{abstract}
In this paper, a
family of infinite dimensional
Lie algebras $\tl$ is introduced and investigated,
called the extended Heisenberg-Virasoro algebra,
 denoted by $\tl$.
These Lie algebras are related to the $N=2$ superconformal algebra
and the Bershadsky-Polyakov algebra.
We study restricted modules and associated vertex algebras
of the Lie algebra $\tl$.
More precisely, we construct its
associated vertex (operator) algebras
$V_{\tl}(\e,0)$,
and show that
the category of vertex algebra $V_{\tl}(\e,0)$-modules
is equivalent to the
category of restricted $\tl$-modules of level $\e$.
 Then we give uniform constructions
of simple restricted $\tl$-modules.
Also, we present several equivalent characterizations of
simple restricted modules over $\tl$.
\end{abstract}

{\bf Keywords} Extended Heisenberg-Virasoro algebra $\cdot$ Restricted module $\cdot$ Irreducible module $\cdot$ Vertex algebra

{\bf Mathematics Subject Classification (2020)} 17B10 $\cdot$ 17B65 $\cdot$ 17B68 $\cdot$ 17B69

\section{Introduction}

Restricted modules (or smooth modules, cf. \cite{KL93})
of Lie algebras are important representations for Lie algebras
and their corresponding vertex algebras.
It is well known that the category of restricted modules
for affine Lie algebras is equivalent to the
category of modules for the corresponding affine vertex algebras (cf. \cite{FZ}, \cite{LL}, etc.).
The study of restricted modules for Lie algebras did not make much
progress until
the work of Mazorchuk and Zhao
in \cite{MZ}.
Since then, (simple) restricted modules for many known Lie algebras are
constructed and characterized
(cf. \cite{CG}, \cite{CHS}, \cite{Ga}, \cite{GG}, \cite{G21}, \cite{GX}, \cite{LPX}, \cite{LPXZ}, etc.).

In this paper, we introduce and study a family of infinite dimensional Lie algebras
which we call extended Heisenberg-Virasoro algebra.
It is a semidirect product of the Heisenberg-Virasoro algebra and its weight modules
of intermediate series.
Heisenberg-Virasoro algebra is introduced and studied in \cite{ACKP},
which is the universal central extension of the semidirect product of the Heisenberg algebra
and the Virasoro algebra.
Representation theory of Heisenberg-Virasoro algebra
has been studied by many authors (cf. \cite{AR}, \cite{Bil},
\cite{CG}, \cite{Ga}, \cite{GW}, \cite{G} for partial references).
Simple restricted modules for Heisenberg-Virasoro algebra are studied in \cite{CG}
and \cite{Ga}.

Let $\CL$ be the Heisenberg-Virasoro algebra which is a Lie algebra spanned by elements $L_{m}, J_{m}$ for $m\in\Z$,
 and central elements
$C_{1}, C_{2}, C_{3}$ (cf. \cite{ACKP}, \cite{Bil}).
For any $\alpha,\beta,F\in\C$, there
is a weight module of intermediate
series $V_{\alpha,\beta,F}=\mbox{span}\{G_{m}\ | \ m\in\Z\}$ for the Lie algebra $\CL$ (cf. \cite{LZ}).
Extended Heisenberg-Virasoro algebra $\L=V_{\alpha,\beta,F}\rtimes\CL$
is a Lie algebra
linearly spanned by
elements $L_{m}, J_{m}, G_{m}$ for $m\in\Z$,
 and central elements
$C_{1}, C_{2}, C_{3}$
with Lie brackets
\begin{equation*}
\begin{aligned}
&[L_{m}, L_{n}]=(m-n)L_{m+n}+\frac{m^{3}-m}{12}\delta_{m+n,0}C_{1},  \ \ \ \ \
[G_{m}, G_{n}]=0,\\
&[L_{m}, J_{n}]=-n J_{m+n}-(m^{2}+m)\delta_{m+n,0}C_{2}, \ \ \ \ \
 [J_{m}, J_{n}]=m\delta_{m+n,0}C_{3},\\
&[L_{m}, G_{n}]=-(\alpha+n+m\beta)G_{m+n},\ \ \ \ \ \ \
[J_{m}, G_{n}]=FG_{m+n}
\end{aligned}
\end{equation*}
for $m,n\in\Z$.
The Lie algebra $\L$ is related to subalgebras of
 $N=2$ superconformal algebra and Bershadsky-Polyakov algebra,
 see Section \ref{sec:2} for the explanation.
For simplicity, denote $\tl=\L$
in the rest
of the introduction.

We first study $\tl$ from the point of view of vertex algebras and their modules,
where we assume $\alpha=\beta$.
Consider the subalgebra
\begin{eqnarray*}
\tl_{+}=\C L_{-1}\oplus\coprod_{i\in\Z_{+}}(\C L_{i}
 \oplus\C J_{i}\oplus\C G_{i})\oplus\sum_{i=1}^{3}\C C_{i}.
\end{eqnarray*}
Let $\C$ be an $\tl_{+}$-module of level $\e$,
form the induced module $V_{\tl}(\e,0)=U(\tl)\otimes_{U(\tl_{+})}\C.$
Then $V_{\tl}(\e,0)$ is a vertex
algebra (cf. \cite{LL}, etc.).
We show that
the category of restricted $\tl$-modules
of level $\e$ is equivalent to the category
of modules for the vertex algebra $V_{\tl}(\e,0)$.

For the study of simple restricted $\tl$-modules,
we divided into two cases.
If $F\neq 0$, for any $d\in\Z_{+}$,
we consider a subalgebra
\begin{eqnarray*}
\tl_d=\sum_{i\in \Z_+}( \C L_i\oplus \C J_i\oplus \C G_{-d+i})\oplus \sum_{i=1}^{3}\C C_i.
\end{eqnarray*}
For any simple $\tl_{d}$-module $V$, we prove that
under certain conditions the induced module
$
\id_{\tl_{d}}^{\tl}(V)=U(\tl)\otimes_{U(\tl_{d})}V
$
is also simple.
If $F=0$ and $C_{3}\neq 0$, for any $\ud=(d_{1},d_{2})\in\Z_{+}^{2}$,
let
\begin{eqnarray*}
\tl_{\ud}=\sum_{i\in\Z_+}(\C L_i\oplus\C G_{-d_1+i}\oplus\C J_{-d_2+i})\oplus\sum_{i=1}^{3}\C C_i.
\end{eqnarray*}
For any simple $\tl_{\ud}$-module $V$
with some conditions,
we show the irreducibility of the induced module
$
\id_{\tl_{\ud}}^{\tl}(V)=U(\tl)\otimes_{U(\tl_{\ud})}V.
$
Then
 we give several equivalent characterizations
of simple restricted $\tl$-modules, as locally nilpotent
or locally finite modules over certain finitely generated subalgebras of $\tl$.
 Also, we show that if
$F\neq 0$ or $F=0,C_{3}\neq0$, all the simple restricted $\tl$-modules
 with some conditions are isomorphic to the induced modules constructed above.
 These give the information of simple weak $V_{\tl}(\e,0)$-modules.
 The automorphism groups and simple weak twisted $V_{\tl}(\e,0)$-modules are studied in \cite{GL}.

  This paper is organized as follows.
In Section \ref{sec:2}, we introduce and study Lie
algebra $\tl$ and its associated vertex algebras $V_{\tl}(\e,0)$.
We prove an equivalence between restricted $\tl$-modules
of level $\e$ and vertex algebra $V_{\tl}(\e,0)$-modules.
In Section \ref{sec:3}, we present constructions of simple $\tl$-modules,
 which show that we can construct
simple restricted modules over $\tl$ from simple modules
over certain subalgebras of $\tl$, see Theorems \ref{idmod1}
and \ref{idmod2}.
In Section \ref{sec:4},
we give several equivalent
characterizations of simple restricted $\tl$-modules.
Moreover, we show that each simple restricted $\tl$-module
satisfying the conditions in Theorems \ref{clasf1} and \ref{clasf2}
is isomorphic to the induced module constructed in Section \ref{sec:3}.

Throughout the paper,
 $\Z$, $\N$, $\Z_+$ and $\C$ are the sets of integers, positive integers,
nonnegative integers and complex numbers, respectively.
The degree $\d(0)$ is not defined and whenever we write $\d(v)$
we mean $v\ne 0$.

\section{Preliminaries and relation with vertex algebras}
 \label{sec:2}
	\def\theequation{2.\arabic{equation}}
	\setcounter{equation}{0}

In this section, we introduce and study the new infinite dimensional
Lie algebra $\tl$ and associated vertex algebras.
We also recall and introduce some notions and total orders for later use.

\vspace{3mm}

{\bf 2.1. Infinite dimensional Lie algebra $\tl$.}
In \cite{ACKP}, the Lie algebra of regular differential operators on $\C^{\times}:=\C\backslash \{0\}$
of degree $\leqslant 1$ were studied from the point of view of cohomology theory.
Its universal central extension $\CL$ is nowadays called the twisted Heisenberg-Virasoro algebra,
which is a Lie algebra spanned by elements $L_{m}, J_{m}$ for $m\in\Z$,
 and central elements
$C_{1}, C_{2}, C_{3}$
with Lie brackets
\begin{equation}\label{eq:HVLie}
\begin{aligned}
&[L_{m}, L_{n}]=(m-n)L_{m+n}+\frac{m^{3}-m}{12}\delta_{m+n,0}C_{1},\\
&[L_{m}, J_{n}]=-n J_{m+n}-(m^{2}+m)\delta_{m+n,0}C_{2}, \\
 &[J_{m}, J_{n}]=m\delta_{m+n,0}C_{3}
\end{aligned}
\end{equation}
for $m,n\in\Z$.
The Lie subalgebra spanned by $\{L_{m},C_{1}\ | \ m\in\Z\}$ is the Virasoro algebra,
and the Lie subalgebra spanned by
$\{J_{m},C_{3}\ | \ m\in\Z\}$ is the
oscillator algebra (also called the Heisenberg algebra).
Note that $\CL$ is the universal central extension
of the semidirect product of the Virasoro algebra
and the oscillator algebra by a 1-dimensional center.

For any $\alpha,\beta,F\in\C$, let $V_{\alpha,\beta,F}=\mbox{span}\{G_{n}\ | \ n\in\Z\}$
be an $\CL$-module with action defined by
\begin{eqnarray}\label{eq:Interm}
\begin{aligned}
&L_{m}G_{n}=-(\alpha+n+m\beta)G_{m+n}, \;\;\;J_{m}G_{n}=FG_{m+n},\\
&C_{1}G_{n}=C_{2}G_{n}=C_{3}G_{n}=0\;\;\;\mbox{for all}\;\;m,n\in\Z.
\end{aligned}
\end{eqnarray}
These are the weight modules of intermediate series over $\CL$ studied in \cite{LZ}.

\bd{}
{\em
For any $\alpha,\beta,F\in\C$,
let $\L=V_{\alpha,\beta,F}\rtimes\CL$
be a semidirect product Lie algebra,
we call it {\em extended Heisenberg-Virasoro algebra}.
Then $\L$ is a Lie algebra linearly spanned by
elements $L_{m}, J_{m}, G_{m}$ for $m\in\Z$,
 and central elements
$C_{1}, C_{2}, C_{3}$
with Lie brackets
\begin{equation}\label{eq:EHVLie}
\begin{aligned}
&[L_{m}, L_{n}]=(m-n)L_{m+n}+\frac{m^{3}-m}{12}\delta_{m+n,0}C_{1},  \ \ \ \ \
[G_{m}, G_{n}]=0,\\
&[L_{m}, J_{n}]=-n J_{m+n}-(m^{2}+m)\delta_{m+n,0}C_{2}, \ \ \ \ \
 [J_{m}, J_{n}]=m\delta_{m+n,0}C_{3},\\
&[L_{m}, G_{n}]=-(\alpha+n+m\beta)G_{m+n},\ \ \ \ \ \ \
[J_{m}, G_{n}]=FG_{m+n}
\end{aligned}
\end{equation}
for $m,n\in\Z$.
}
\ed
Clearly, if $F=0$, then the element $J_{0}$ is also in the center.
If $\alpha=\beta=F=0$,
then $G_{0}$ is also a central element.

We notice the following Lie algebra isomorphisms.
\bp{aut-L}
Let $\alpha,\beta,F\in\C$.
\vspace{3mm}
\\
{\rm(1)}
$\CL_{\alpha,\beta,F}\cong \CL_{\alpha+p,\beta,F}$ for any $p\in\Z$;\\
{\rm (2)} Suppose $\alpha\notin\Z$, then
$\CL_{\alpha,0,0}\cong \CL_{\alpha,1,0}$;\\
{\rm(3)} For $F\in\C^{\times}$,
we have
$\CL_{\alpha,\beta,F}\cong \CL_{\alpha,\beta,1}$;\\
{\rm(4)} If $C_{3}=2C_{2}$, then $\CL_{0,0,1}\cong\CL_{0,-1,-1}$.
\ep
\begin{proof}
For any $p\in\Z$,
let $\vf_{1}$ be a linear map from
$\CL_{\alpha,\beta,F}$ to $\CL_{\alpha+p,\beta,F}$ defined by
$\vf_{1}(C_{i})=C_{i}$ for $i=1,2,3$,
 and for $m\in\Z$,
\begin{eqnarray}
\vf_{1}:
L_{m}\mapsto L_{m},\;\;J_{m}\mapsto J_{m},\;\;G_{m}\mapsto G_{m-p}.
\end{eqnarray}
Suppose $\alpha\notin\Z$.
Let $\vf_{2}$ be a linear map
from $\CL_{\alpha,0,0}$
to $\CL_{\alpha,1,0}$
defined by
$\vf_{2}(C_{i})=C_{i}$
for $i=1,2,3$,
 and for $m\in\Z$,
\begin{eqnarray}\vf_{2}:
L_{m}\mapsto L_{m},\;\;J_{m}\mapsto J_{m},\;\;G_{m}\mapsto (\alpha+m)G_{m}.
\end{eqnarray}
For any $F\in\C^{\times}$,
let $\vf_{3}$ be a linear map
from $\CL_{\alpha,\beta,F}$
to $\CL_{\alpha,\beta,1}$
defined by
\begin{eqnarray}
\begin{aligned}
\vf_{3}:\;\;
&C_{1}\mapsto C_{1},\;\;C_{2}\mapsto FC_{2},\;\;C_{3}\mapsto F^{2}C_{3},\\
&L_{m}\mapsto L_{m},\;\;J_{m}\mapsto FJ_{m},\;\;G_{m}\mapsto G_{m}
\;\;\mbox{for}\;m\in\Z.
\end{aligned}
\end{eqnarray}
If $C_{3}=2C_{2}$,
let $\vf_{4}$ be a linear map
from $\CL_{0,0,1}$
to $\CL_{0,-1,-1}$
defined by
$\vf_{4}(C_{i})=C_{i}$
for $i=1,2$,
 and for $m\in\Z$,
\begin{eqnarray}\vf_{4}:
L_{m}\mapsto L_{m}+mJ_{m},\;\;J_{m}\mapsto -J_{m}+2\delta_{m,0}C_{2},
\;\;G_{m}\mapsto -G_{m}.
\end{eqnarray}
It is straightforward to check that $\vf_{i}$,
$i=1,2,3,4$,
are Lie algebra isomorphisms.
\end{proof}

The Lie algebra $\L$ is closely related to $N=2$ superconformal algebra
and Bershadsky-Polyakov algebra.

\vspace{3mm}

{\bf 2.2. Relation with $N=2$ superconformal algebra.}
$N=2$ superconformal algebra is an infinite-dimensional Lie superalgebra,
related to supersymmetry, that occurs in string theory and conformal
field theory as a gauge algebra of the $U(1)$ fermionic string (cf. \cite{A}, \cite{K}, etc.) .
It has important applications in mirror symmetry.
Recall that the (Ramond) $N=2$ superconformal
algebra is a Lie superalgebra with basis of even elements
$L_{m},J_{m}$ for $m\in\Z$,
central element ${\bf c}$,
odd elements $G_{m}^{\pm}$ for $m\in\Z$,
and relations
\begin{eqnarray}\label{eq:N=2Lie}
\begin{aligned}
&[L_{m}, L_{n}]=(m-n)L_{m+n}+\frac{m^{3}-m}{12}\delta_{m+n,0}{\bf c},  \ \ \ \ \
\{G^{\pm}_{m}, G^{\pm}_{n}\}=0,\\
&[L_{m}, J_{n}]=-n J_{m+n}, \ \ \ \ \
 [J_{m}, J_{n}]=m\delta_{m+n,0}\frac{{\bf c}}{3},\\
&[L_{m}, G^{\pm}_{n}]=(\frac{m}{2}-n)G^{\pm}_{m+n},\ \ \ \ \ \ \
[J_{m}, G^{\pm}_{n}]=\pm G^{\pm}_{m+n},\\
&\{G_{m}^{+},G_{n}^{-}\}=L_{m+n}+\frac{1}{2}(m-n)J_{m+n}+\frac{{\bf c}}{6}(r^{2}-\frac{1}{4})\delta_{m+n,0}
\end{aligned}
\end{eqnarray}
for $m,n\in\Z$.
It is easy to see that if we take $C_{1}={\bf c}$, $C_{3}=\frac{1}{3}{\bf c}=\frac{1}{3}C_{1}$ and $C_{2}=0$, then
Lie algebras $\CL_{0,-\frac{1}{2},\pm 1}$
are closely related to the subalgebras
(generated by $L,J,G^{+}$, or $L,J,G^{-}$) of the (Ramond) $N=2$ superconformal
algebra.

\vspace{3mm}

{\bf 2.3. Relation with Bershadsky-Polyakov algebra.}
Bershadsky-Polyakov algebra is a vertex algebra freely generated
by four fields $L(x),J(x),G^{\pm}(x)$
and satisfying certain OPEs (cf. \cite{B}, \cite{P}, etc.),
where $x$ is a formal variable.
We recall from \cite{Ara13} the OPEs of the Bershadsky-Polyakov algebra.
Expand the four fields as $L(x)=\sum_{n\in\Z}L_{n}x^{-n-2}$
and
\begin{eqnarray}
J(x)=\sum_{n\in\Z}J_{n}x^{-n-1},
G^{+}(x)=\sum_{n\in\Z}G^{+}_{n}x^{-n-1},
G^{-}(x)=\sum_{n\in\Z}G^{-}_{n}x^{-n-2}.
\end{eqnarray}
Then
\begin{eqnarray}\label{eq:BPLie}
\begin{aligned}
&[L_{m}, L_{n}]=(m-n)L_{m+n}+\frac{m^{3}-m}{12}\delta_{m+n,0}c(k),  \ \ \ \ \
[G^{\pm}_{m}, G^{\pm}_{n}]=0,\\
&[L_{m}, J_{n}]=-n J_{m+n}-(m^{2}+m)\delta_{m+n,0}\frac{2k+3}{6}, \
 [J_{m}, J_{n}]=m\delta_{m+n,0}\frac{2k+3}{3},\\
&[L_{m}, G^{+}_{n}]=-nG^{+}_{m+n},\ \ \ \ \ \ \
[J_{m}, G^{+}_{n}]=G^{+}_{m+n},\\
&[L_{m}, G^{-}_{n}]=(m-n)G^{-}_{m+n},\ \ \ \ \ \ \
[J_{m}, G^{-}_{n}]=-G^{-}_{m+n}
\end{aligned}
\end{eqnarray}
for $m,n\in\Z$, where $c(k)=-\frac{4(k+1)(2k+3)}{k+3}$, $k\neq-3,k\in\C$.
For the relation between $G_{m}^{+}$ and $G_{n}^{-}$ (we don't need it here), see \cite{Ara13} for example.
For the extended Heisenberg-Virasoro algebra $\L$,
if we take
$C_{1}=c(k)$,
$C_{2}=\frac{2k+3}{6}$,
$C_{3}=2C_{2}=\frac{2k+3}{3}$
for some $-3\neq k\in\C$,
then Lie algebras $\CL_{0,0,1}$
and $\CL_{0,-1,-1}$
are the same as Lie algebras spanned by $L,J,G^{+}$ and $L,J,G^{-}$ from the Bershadsky-Polyakov algebra.

In the following, for simplicity, we
denote $\tl=\L$.
For future convenience, we recall and introduce some notions and total orders.

\vspace{3mm}

{\bf 2.4. Notions and total orders.}
Let $W$ be a general vector space.
  Set
  $$
\mathcal{E}(W)=\mbox{Hom}(W , W((x)))\subset(\mbox{EndW})[[x,x^{-1}]],
$$
where
$$W((x))=\{\sum\limits_{n\in\Z}w_{n}x^{n}\ | \ w_{n}\in W,
w_{n}=0\;\mbox{for}\;n\;\mbox{sufficiently negative}\}.$$
The identity operator on $W$, denoted by $\textbf{1}_{W}$, is a special
element of $\mathcal{E}(W).$
Recall that $L(x)=\sum_{n\in\Z}L_{n}x^{-n-2}$,
$J(x)=\sum_{n\in\Z}J_{n}x^{-n-1}$.
Denote $G(x)=\sum\limits_{n\in\Z}G_{n}x^{-n-1}$.

 \bd{restricted and level}
{\em An $\tl$-module $W$ is called {\em restricted} if for any $w\in W$,
$L_{m}w=J_{m}w=G_{m}w=0$ for $m$ sufficiently large,
or equivalently, if $L(x), J(x), G(x)\in {\mathcal{E}}(W)$.
We say an $\tl$-module $W$ is of {\em level} $\e$
if the central element $C_{i}$ acts as
a scalar $\ell_{i}$ for $i=1,2,3.$
}
\ed

\bd{}
{\em Let $L$ be a Lie algebra and $W$ be an $L$-module.
 For $x \in L$, we say that $x$ acts {\em locally nilpotently}
on $W$ if for any $w\in W$ there exists $n\in \N$ such that $x^{n}w=0$
and that $x$ acts {\em locally finitely} on $W$
if for any $w\in W$ we have $\mbox{dim}(\sum_{n\in \N}{\C x^nw})<+\infty$.
 We say that $L$ acts {\em locally nilpotently} on $W$ if for any $w\in W$
  there exists $n\in \N$ such that $L^{n}w=0$
and $L$ acts {\em locally finitely} on $W$ if for any $w\in W$
 we have $\mbox{dim}(\sum_{n\in \N}{L^nw})<+\infty$.
}\ed

Let $\M$ be the set of all infinite vectors of the form
${\bf i}=(...,i_2,i_1)$ with entries in $\Z_+$
such that the number of nonzero entries is finite.
Let {\bf 0} denote the element $(...,0,0)\in \M$ and
$\es_i$ denote the element $(...,0,1,0,...,0)\in \M$ for $i\in \N$,
 where 1 is in the $i$-th position from right.
 For any ${\bf i}\in \M$, denote
\begin{equation}
  {\bf w(i)}=\sum_{s\in \N}{s\cdot i_s},\ \ \ \ \ \ \
  {\bf d(i)}=\sum_{s\in \N}{i_s},
\end{equation}
which are nonnegative integers.

Let $>$ and $\succ$ be the {\em lexicographical} and
 {\em reverse lexicographical} total orders on $\M$ respectively,
 that is, for any ${\bf i},{\bf j} \in \M$,
\begin{equation*}
{\bf j}>{\bf i}\ \Leftrightarrow\ {\rm there\ exists}\ r\in \N\ {\rm such\ that}\
(j_s=i_s ,\forall\; s>r)\ {\rm and}\  j_r>i_r,
\end{equation*}
and
\begin{equation*}
{\bf j}\succ {\bf i}\ \Leftrightarrow \ {\rm there\ exists}\ r\in \N\ {\rm such\ that}\
 (j_s=i_s, \forall\; 1\leqslant s<r)\ {\rm and}\  j_r>i_r.
\end{equation*}

\vspace{3mm}

{\bf 2.5 Vertex algebras associated to $\tl$.}
We now construct vertex algebras $V_{\tl}(\e,0)$
associated to the Lie algebra $\tl=\L$,
for which we require that $\alpha=\beta$.
For the case of $\alpha\neq\beta$,
by rewriting basis of $\tl$, see Proposition \ref{aut-L},
perhaps it is also possible to construct the corresponding vertex algebras,
which we don't consider it here.
For the notions and notations about vertex (operator) algebras
and their modules, we follow \cite{LL}.
We show that the category of restricted $\tl$-modules of level $\e$
 is equivalent to the
category of modules of the vertex algebra $V_{\tl}(\e,0)$.
In particular, simple modules
for the vertex algebra $V_{\tl}(\e,0)$
correspond to simple restricted $\tl$-modules
of level $\e$.

It is clear that $\tl$ is a $\Z$-graded Lie algebra
$\tl=\coprod_{n\in\Z}\tl_{(n)}$ with
\begin{eqnarray}
\begin{aligned}
&\tl_{(0)}=\C L_{0}\oplus \C J_{0}\oplus \C G_{0}\oplus\sum_{i=1}^{3}\C C_{i},\\
&\tl_{(n)}=\C L_{-n}\oplus \C J_{-n}\oplus\C G_{-n}\;\;\mbox{for}\;\;n\neq 0.
\end{aligned}
\end{eqnarray}
Because
$[L_{0}, G_{n}]=-(\alpha+n)G_{n},$
 this grading is not given by $\ad L_{0}$-eigenvalues
unless $\alpha=0$.
Let
 \begin{eqnarray}
 \begin{aligned}
 &\tl_{+}=\C L_{-1}\oplus\coprod_{n\in\Z_{+}}(\C L_{n}
 \oplus\C J_{n}\oplus\C G_{n})\oplus\sum_{i=1}^{3}\C C_{i},\\
& \tl_{-}=\coprod_{n\leqslant -2}(\C L_{n}
 \oplus\C J_{n}\oplus\C G_{n})\oplus\C J_{-1}\oplus\C G_{-1}.
 \end{aligned}
 \end{eqnarray}
 Note that
 $[L_{-1},G_{0}]=-(\alpha-\beta)G_{-1}.$
 Since we assume $\alpha=\beta$, $\tl_{+}$ and $\tl_{-}$
 are subalgebras of $\tl$, and $\tl=\tl_{+}\oplus\tl_{-}$.

 Let $\ell_{i}$, $i=1,2, 3$ be any complex numbers.
Let $\C$ be an $\tl_{+}$-module with $C_{i}$ acting
as the scalar $\ell_{i}$ for $i=1,2, 3$, and
all other elements
acting trivially.
Form the induced module
\begin{eqnarray} \label{eq:ind}
V_{\tl}(\e,0)=U(\tl)\otimes_{U(\tl_{+})}\C.
\end{eqnarray}
Identify $\C$ as a subspace of $V_{\tl}(\e,0)$
and set ${\bf 1} =1\in\C\subset V_{\tl}(\e,0)$.
 Denote $\omega=L_{-2}{\bf 1}$, $J=J_{-1}{\bf 1}$ and $G=G_{-1}{\bf 1}.$
Then
$V_{\tl}(\e,0)=\coprod_{n\geqslant 0}V_{\tl}(\e,0)_{(n)}$,
where $V_{\tl}(\e,0)_{(0)}=\C$,
 and $V_{\tl}(\e,0)_{(n)}$, $n\geqslant 1$, has a basis consisting of vectors
$$G_{-n_{1}}\cdots G_{-n_{p}}J_{-k_{1}}\cdots J_{-k_{s}}L_{-m_{1}}\cdots L_{-m_{r}}{\bf{1}} $$
for
$r, s,p\geqslant 0$,
$m_{1}\geqslant\cdots\geqslant m_{r}\geqslant 2$, $k_{1}\geqslant\cdots\geqslant k_{s}\geqslant 1$,
 $n_{1}\geqslant\cdots\geqslant n_{p}\geqslant 1$ with
 $\sum\limits_{i=1}^{r}m_{i}+\sum\limits_{j=1}^{s}k_{j}+\sum\limits_{t=1}^{p}n_{t}=n.$

 The proof of the following result is standard
 (cf. Theorem 5.7.1 of \cite{LL}, etc.).
\bt{VA1}
Let $\ell_{1},\ell_{2},\ell_{3}\in\C$ be arbitrary.
Suppose that $\alpha=\beta$.
\begin{itemize}
\item [(1)]
 $V_{\tl}(\e,0)$ is a vertex algebra,
which is uniquely determined by the condition that ${\bf 1}$ is the vacuum vector
 and
$Y(L_{-2}{\bf 1},x)=L(x),\;\;
Y(J_{-1}{\bf 1},x)=J(x),\;\;
Y(G_{-1}{\bf 1},x)=G(x)$.
\item [(2)]
If $\alpha=\beta=0$,
then $V_{\tl}(\e,0)$ is a vertex operator algebra with
$\omega=L_{-2}{\bf 1}$ the conformal vector.
\end{itemize}
\et

As an $\tl$-module,
$V_{\tl}(\e,0)$ is generated by ${\bf 1}$ with the
relations $C_{i}=\ell_{i}$ and
$L_{n}{\bf 1}=J_{m}{\bf 1}=G_{m}{\bf 1}=0$ for $n\geqslant -1, m\geqslant 0,$ $i=1,2, 3$.
In fact, $V_{\tl}(\e,0)$ is {\em universal}
in the sense that for any $\tl$-module $W$
of level $\e$ equipped with a vector $w\in W$
such that $L_{n}w=J_{m}w=G_{m}w=0$ for $n\geqslant -1, m\geqslant 0,$
there exists a unique
$\tl$-module homomorphism from $V_{\tl}(\e,0)$
to $W$ sending ${\bf 1}$ to $w$.

We have the following theorem.
\bt{Mods}
There is a one-to-one correspondence
between restricted $\tl$-modules $W$
of level $\e$
and
modules of the vertex algebra
$V_{\tl}(\e,0)$ with
$
Y_{W}(L_{-2}{\bf 1},x)=L(x)$,
$Y_{W}(J_{-1}{\bf 1},x)=J(x)$,
$Y_{W}(G_{-1}{\bf 1},x)=G(x).
$
\et
\begin{proof}
Suppose that $W$ is a restricted $\tl$-module
of level $\e.$
Consider the set
$
U_{W}=\{L(x) ,J(x), G(x), {\bf 1}_{W} \}.
$
Then $U_{W}$ is a subset of $\mathcal{E}(W)$.
From the defining relations (\ref{eq:EHVLie}) of $\tl$
we see that $U_{W}$ is a local subset.
By Theorem 5.5.18 of \cite{LL},
$U_{W}$ generates a vertex algebra $\langle U_{W}\rangle$ with $W$ a natural faithful module.
Moreover, $\langle U_{W}\rangle$ is an $\tl$-module of level $\e$,
with $ L_{n}, J_{n}, G_{n}$ acting as
$L(x)_{n+1}, J(x)_{n}, G(x)_{n}$
for $n\in\Z$
and with the relations $L(x)_{m}{\bf 1}_{W}=J(x)_{n}{\bf 1}_{W}=G(x)_{n}{\bf 1}_{W}=0$
for $m\geq -1$, $n\geq 0$.
By the universal property of $V_{\tl}(\e,0)$,
there exists a unique $\tl$-module
homomorphism $\psi$ from
$V_{\tl}(\e,0)$ to $\langle U_{W}\rangle$
such that $\psi({\bf 1})={\bf 1}_{W}$.
Then we have
\begin{eqnarray*}
\psi(\omega_{n}v)= L(x)_{n}\psi(v),\;\; \psi(J_{n}v)= J(x)_{n}\psi(v)
\;\;\mbox{and}\;\;\psi(G_{n}v)= G(x)_{n}\psi(v)
\end{eqnarray*}
for
 $n\in\Z, v\in V_{\tl}(\e,0)$.
The existence and uniqueness of $V_{\tl}(\e,0)$-module
 structure on $W$ now immediately
follows from Theorem 5.7.6 of \cite{LL}.

Conversely, let $W$ be a module for the vertex algebra $V_{\tl}(\e,0)$.
For $i\geq 0$, we have
\begin{eqnarray}
\begin{aligned}
& (L_{-2}{\bf 1})_{i}L_{-2}{\bf 1}
 = (i+1)L_{i-3}
{\bf 1}+\delta_{i-3,0}\frac{(i-1)^{3}-(i-1)}{12}C_{1}{\bf 1},\\
& (L_{-2}{\bf 1})_{i}J_{-1}{\bf 1}
= J_{i-2}
{\bf 1}-\delta_{i-2,0}((i-1)^{2}+(i-1))C_{2}{\bf 1},\\
&(L_{-2}{\bf 1})_{i}G_{-1}{\bf 1}=-(\alpha-1+(i-1)\beta)G_{i-2}{\bf 1},\;\;\;
\;(G_{-1}{\bf 1})_{i}G_{-1}{\bf 1}=0,\\
&(J_{-1}{\bf 1})_{i}J_{-1}{\bf 1}= i\delta_{i-1,0}C_{3}{\bf 1},\;\;\;\;\;\;
(J_{-1}{\bf 1})_{i}G_{-1}{\bf 1}=FG_{i-1}{\bf 1}.
\end{aligned}
\end{eqnarray}
Note that (cf. Lemma 2.3.5 of \cite{L2})
$$[Y_{W}(a,x_{1}),Y_{W}(b,x_{2})]=\sum_{i\geq 0}\frac{1}{i!}Y_{W}(a_{i}b,x_{2})
\left(\frac{\partial}{\partial x_{2}}\right)^{i}x_{1}^{-1}\delta\left(\frac{x_{2}}{x_{1}}\right).$$
By a direct calculation, we show that $W$ is an $\tl$-module of {\em level} $\e$
with
$L(x)=Y_{W}(L_{-2}{\bf 1},x)$, $J(x)=Y_{W}(J_{-1}{\bf 1},x)$
and $G(x)=Y_{W}(G_{-1}{\bf 1},x)$.
Since $W$ is a $V_{\tl}(\e,0)$-module,
by definition, $Y_{W}(L_{-2}{\bf 1},x)$,
$Y_{W}(J_{-1}{\bf 1},x)$, $Y_{W}(G_{-1}{\bf 1},x)\in\mathcal{E}(W)$.
Therefore, $W$ is a restricted $\tl$-module of level $\e$.
\end{proof}

In the rest of this paper, we study simple restricted $\tl$-modules
in the general setting, where we do not require $\alpha=\beta$.

\section{Construction of simple $\tl$-modules}
\label{sec:3}
\def\theequation{3.\arabic{equation}}
	\setcounter{equation}{0}

In this section, we present two constructions
of simple $\tl$-modules depending on whether $F$
is zero or not.
These modules are obtained from simple modules
over certain subalgebras of $\tl$.

\vspace{3mm}

{\bf 3.1. The case of $F\neq 0$.}

\vspace{3mm}

Let $d\in \Z_+$ be arbitrary.
Consider the subalgebra
\begin{equation}\label{eq:tld}
\tl_d=\sum_{i\in \Z_+}( \C L_i\oplus \C J_i\oplus \C G_{-d+i})\oplus \sum_{i=1}^{3}\C C_i.
\end{equation}
For any simple $\tl_{d}$-module $V$,
we have the induced $\tl$-module
\begin{equation*}
\id_{\tl_{d}}^{\tl}(V) =U(\tl)\otimes _{U(\tl_d)}V.
\end{equation*}
For ${\bf k,i,j} \in \M$, denote
\begin{equation}\nonumber
L^{\bf k}G^{\bf i}J^{\bf j}
=\cdots L_{-2}^{k_2}L_{-1}^{k_1}\cdots G_{-d-2}^{i_2}G_{-d-1}^{i_1}\cdots J_{-2}^{j_2}J_{-1}^{j_1}\in U(\tl).
\end{equation}
By the PBW Theorem, each element of $\id_{\tl_{d}}^{\tl}(V)$ can be uniquely written in the form
\begin{equation}\label{eq:idelt}
 \sum_{{\bf k,i,j}\in \M}{L^{\bf k}G^{\bf i}J^{\bf j}v_{\bf k,i,j}},
\end{equation}
where all $v_{\bf k,i,j} \in V$ and only finitely many of them are nonzero.

\bd{pto}
{\em Define a {\em principal total order} $({\bf o})$ on $\M\times\M\times\M$,
 still denoted by $\succ$, as follows:
 for any ${\bf k,i,j,a,b,c} \in \M $},
\begin{equation}\nonumber
({\bf k,i,j}) \succ ({\bf a,b,c})\Leftrightarrow
  ({\bf k,i,j,w(k),w(k+i+j)}) \succ ({\bf a,b,c,w(a),w(a+b+c)}).
\end{equation}
\ed

For any $v\in\id_{\tl_{d}}^{\tl}(V)$ written in the form of (\ref{eq:idelt}),
 we denote by $\mbox{supp}(v)$ the set of all $({\bf k,i,j)}\in \M\times\M\times\M$
 such that $v_{\bf k,i,j}\ne 0$.
 For a nonzero $v\in \id_{\tl_{d}}^{\tl}(V)$,
  let $\d(v)$ denote the maximal
   (with respect to the principal total order $({\bf o})$ on $\M\times\M\times\M$)
   element in $\mbox{supp}(v)$, called the {\em degree} of $v$.

 \bt{idmod1}
Suppose that $F\ne 0$.
Let $d\in \Z_+$ and $V$ be a simple $\tl_d$-module.
Assume that there exist $k,l\in\Z_{+}$ with $l>k+d>0$ such that
the following conditions hold.
\vspace{3mm}\\
(a) $G_{k}$ and $L_{l}$ act injectively on $V$,\\
(b) $G_{k+i}V=J_{k+d+i}V=L_{l+i}V=0$ for all $i\in\N$.
\vspace{2mm}\\
Then ${\rm Ind}_{\tl_{d}}^{\tl}(V)$ is a simple $\tl$-module.
\et

Theorem \ref{idmod1} is obtained by using the following lemmas
\ref{lemred1.1}, \ref{lemred1.2} and \ref{lemred1.3} repeatedly.
\bl{lemred1.1}
Let $d\in\Z _+$ and $V$ be an $\tl_d$-module
 satisfying the conditions in Theorem \ref{idmod1}.
 For any $v\in {\rm Ind}_{\tl_{d}}^{\tl}(V)\backslash V$,
 we write $v$ in the form of (\ref{eq:idelt}).
 Suppose ${\rm deg}(v)= ({\bf k,i,j})$ with ${\bf j\ne 0}$.
 Let $r={\rm min}\{s:j_s \ne 0\}>0$,
then ${\rm deg}(G_{k+r}v) =({\bf k,i,j}-\es_{r})$.
\el
\begin{proof}
It suffices to consider those $v_{\bf k^{\p},i^{\p},j^{\p}}$ with
 $G_{k+r}L^{\bf k^{\p}}G^{\bf i^{\p}}J^{\bf j^{\p}}v_{\bf k^{\p},i^{\p},j^{\p}}\ne 0$.
 Noticing that $G_{k+r}v_{\bf k^{\p},i^{\p},j^{\p}}=0$
  for any such $({\bf k^{\p},i^{\p},j^{\p}}) \in \mbox{supp}(v)$.
  We have
\begin{equation}\nonumber
G_{k+r}L^{\bf k^{\p}}G^{\bf i^{\p}}J^{\bf j^{\p}}v_{\bf k^{\p},i^{\p},j^{\p}}
=[G_{k+r},L^{\bf k^{\p}}]G^{\bf i^{\p}}J^{\bf j^{\p}}v_{\bf k^{\p},i^{\p},j^{\p}}
\;+\;
L^{\bf k^{\p}}G^{\bf i^{\p}}[G_{k+r},J^{\bf j^{\p}}]v_{\bf k^{\p},i^{\p},j^{\p}}.
\end{equation}

\vspace{3mm}

{\bf Claim \ref{lemred1.1}.1:}
If ${\bf j^{\p}\neq 0}$,
then
${\rm deg}(L^{\bf k^{\p}}G^{\bf i^{\p}}[G_{k+r},J^{\bf j^{\p}}]v_{\bf k^{\p},i^{\p},j^{\p}})
\preceq({\bf k,i,j}-\es_{r})$,
where the equality holds if and only if ${\bf k^{\p}=k,i^{\p}=i,j^{\p}=j}$.

{\em Proof of Claim \ref{lemred1.1}.1:}  We have
\begin{eqnarray}
&L^{\bf k^{\p}}G^{\bf i^{\p}}[G_{k+r},J^{\bf j^{\p}}]v_{\bf k^{\p},i^{\p},j^{\p}} \nonumber\\
&=\sum\limits_{{\j^{*}\in\M}, y\geqslant 1}
L^{\bf k^{\p}}G^{\bf i^{\p}}J^{\j^{*}}[G_{k+r},J_{-y}]J^{{\j^{\p}}-{\j^{*}}-\es_{y}}v_{\bf k^{\p},i^{\p},j^{\p}}
\nonumber\\
&=-F\sum\limits_{{\j^{*}\in\M}, y\geqslant 1}
L^{\bf k^{\p}}G^{\bf i^{\p}}J^{\j^{*}}G_{k+r-y}J^{{\j^{\p}}-{\j^{*}}-\es_{y}}v_{\bf k^{\p},i^{\p},j^{\p}}
\end{eqnarray}
Note that $F\neq 0$.
Suppose $G_{k+r}$ makes commutator with $J_{-y}$, denote
\begin{eqnarray}
\d(L^{\bf k^{\p}}G^{\bf i^{\p}}J^{\j^{*}}G_{k+r-y}J^{{\j^{\p}}-{\j^{*}}-\es_{y}}v_{\bf k^{\p},i^{\p},j^{\p}})
=({\k}_{y},{\i}_{y},{\j}_{y}).
\end{eqnarray}

We now determine $({\k}_{y},{\i}_{y},{\j}_{y})$ according to $y$.
\begin{itemize}
\item[(1)] For $y> k+d+r$, we have
$({\k}_{y},{\i}_{y},{\j}_{y})=({\k^{\p}},{\i^{\p}}+\es_{y-r-k-d},{\j^{\p}}-\es_{y})$.
So we get (note that $k+d>0$)
\begin{eqnarray*}
{\bf w}({\k}_{y}+{\i}_{y}+{\j}_{y})
={\bf w(\k^{\p}+\i^{\p}+\j^{\p})}-k-d-r
<{\bf w(\k+\i+\j)}-r.
\end{eqnarray*}
\item[(2)]For $r<y\leqslant k+d+r$,
it is straightforward to check that
$({\k}_{y},{\i}_{y},{\j}_{y})\preceq({\k^{\p}},{\i^{\p}},{\j^{\p}}-\es_{y})$.
Hence
${\bf w}({\k}_{y}+{\i}_{y}+{\j}_{y})
\leqslant {\bf w(\k^{\p}+\i^{\p}+\j^{\p})}-y
<{\bf w(\k+\i+\j)}-r $.
\item[(3)]For $1\leqslant y<r$,
we have
${\bf w}({\k}_{y}+{\i}_{y}+{\j}_{y})\leqslant{\bf w(\k^{\p}+\i^{\p}+\j^{\p})}-r$.
If ${\bf w(\k^{\p}+\i^{\p}+\j^{\p})}<{\bf w(\k+\i+\j)}$, then
 $({\k}_{y},{\i}_{y},{\j}_{y})\prec({\bf k,i,j}-\es_{r})$.
If ${\bf w(\k^{\p}+\i^{\p}+\j^{\p})}={\bf w(\k+\i+\j)}$,
we have ${\bf w}({\k^{\p}})<{\bf w}(\k)$.
So ${\bf w}({\k}_{y})={\bf w}({\k^{\p}})<{\bf w}(\k)$.
Therefore, $({\k}_{y},{\i}_{y},{\j}_{y})\prec({\bf k,i,j}-\es_{r})$.
\item[(4)] For $y=r$, we have
$({\k}_{y},{\i}_{y},{\j}_{y})=({\k^{\p}},{\i^{\p}},{\j^{\p}}-\es_{r})$.
If ${\bf w(\k^{\p}+\i^{\p}+j^{\p})}<{\bf w(\k+\i+\j)}$, then
${\bf w}({\k}_{y}+{\i}_{y}+{\j}_{y})
={\bf w(\k^{\p}+\i^{\p}+\j^{\p})}-r<{\bf w}({\bf k+i+j}-\es_{r})$.
If ${\bf w(\k^{\p}+\i^{\p}+j^{\p})}={\bf w(\k+\i+\j)}$
and ${\bf w}(\k^{\p})<{\bf w}(\k)$, we have
${\bf w}({\k}_{y}+{\i}_{y}+{\j}_{y})
={\bf w(\k^{\p}+\i^{\p}+\j^{\p})}-r={\bf w}({\bf k+i+j}-\es_{r})$,
but ${\bf w}({\k}_{y})={\bf w(\k^{\p})}<{\bf w}(\k)$.
So we always have
$({\k}_{y},{\i}_{y},{\j}_{y})\prec({\bf k,i,j}-\es_{r})$.
If ${\bf w(\k^{\p}+\i^{\p}+j^{\p})}={\bf w(\k+\i+\j)}$
and ${\bf w}(\k^{\p})={\bf w}(\k)$,
then we must have $\j^{\p}\preceq\j$.
Therefore,
$({\k}_{y},{\i}_{y},{\j}_{y})\preceq({\bf k,i,j}-\es_{r})$,
and the equality holds if and only if ${\bf k^{\p}=k,i^{\p}=i,j^{\p}=j}$.
\end{itemize}

\vspace{3mm}

{\bf Claim \ref{lemred1.1}.2:}
If ${\bf k^{\p}\neq 0}$, then
${\rm deg}([G_{k+r},L^{\bf k^{\p}}]G^{\bf i^{\p}}J^{\bf j^{\p}}v_{\bf k^{\p},i^{\p},j^{\p}})
\prec({\bf k,i,j}-\es_{r})$.

{\em Proof of Claim \ref{lemred1.1}.2:}
Note that
\begin{eqnarray*}
&[G_{k+r},L^{\bf k^{\p}}]G^{\bf i^{\p}}J^{\bf j^{\p}}v_{\bf k^{\p},i^{\p},j^{\p}}\\
&=\sum\limits_{{\bf k^{*}}\in\M,x\geqslant 1}
L^{\bf k^{*}}[G_{k+r},L_{-x}]L^{{\bf k^{\p}-k^{*}}-\es_{x}}G^{\bf i^{\p}}J^{\bf j^{\p}}v_{\bf k^{\p},i^{\p},j^{\p}}\\
&=\sum\limits_{{\bf k^{*}}\in\M,x\geqslant 1}(\alpha+k+r-x\beta)
L^{\bf k^{*}}G_{k+r-x}L^{{\bf k^{\p}-k^{*}}-\es_{x}}G^{\bf i^{\p}}J^{\bf j^{\p}}v_{\bf k^{\p},i^{\p},j^{\p}}.
\end{eqnarray*}
Suppose $G_{k+r}$ makes commutator with $L_{-x}$, denote
$${\rm deg}(L^{\bf k^{*}}G_{k+r-x}L^{{\bf k^{\p}-k^{*}}-\es_{x}}G^{\bf i^{\p}}J^{\bf j^{\p}}v_{\bf k^{\p},i^{\p},j^{\p}})
=({\bf k}_{x},{\bf i}_{x},{\bf j}_{x}).$$

\begin{itemize}
\item[(i)] For $x> k+d+r$, we have
$({\k}_{x},{\i}_{x},{\j}_{x})=({\k^{\p}}-\es_{x},{\i^{\p}}+\es_{x-r-k-d},\j^{\p})$.
Then (note that $k+d>0$)
\begin{eqnarray*}
{\bf w}({\k}_{x}+{\i}_{x}+{\j}_{x})
={\bf w(\k^{\p}+\i^{\p}+\j^{\p})}-k-d-r
<{\bf w(\k+\i+\j)}-r.
\end{eqnarray*}
\item[(ii)] For $r<x\leqslant k+d+r$,
we have
$({\k}_{x},{\i}_{x},{\j}_{x})\preceq({\k^{\p}}-\es_{x},{\i^{\p}},\j^{\p})$.
Hence,
${\bf w}({\k}_{x}+{\i}_{x}+{\j}_{x})
\leqslant {\bf w(\k^{\p}+\i^{\p}+\j^{\p})}-x
<{\bf w(\k+\i+\j)}-r
={\bf w}({\bf k+i+j}-\es_{r})$.
\item[(iii)]
 For $x=r$, we have
$({\k}_{x},{\i}_{x},{\j}_{x})=({\k^{\p}}-\es_{r},{\i^{\p}},\j^{\p})$.
If ${\bf w(\k^{\p}+\i^{\p}+j^{\p})}<{\bf w(\k+\i+\j)}$, then
${\bf w}({\k}_{x}+{\i}_{x}+{\j}_{x})
={\bf w(\k^{\p}+\i^{\p}+\j^{\p})}-r<{\bf w}({\bf k,i,j}-\es_{r})$.
If ${\bf w(\k^{\p}+\i^{\p}+j^{\p})}={\bf w(\k+\i+\j)}$, we get
${\bf w}({\k}_{x}+{\i}_{x}+{\j}_{x})
={\bf w(\k^{\p}+\i^{\p}+\j^{\p})}-r={\bf w}({\bf k+i+j}-\es_{r})$,
but ${\bf w} ({\k}_{x})={\bf w(\k^{\p})}-r<{\bf w}(\k)$.
So we always have
$({\k}_{x},{\i}_{x},{\j}_{x})\prec({\bf k,i,j}-\es_{r})$.
\item[(iv)] For $1\leqslant x<r$, we have
\begin{eqnarray*}
&L^{\bf k^{*}}[G_{k+r},L_{-x}]L^{{\bf k^{\p}-k^{*}}-\es_{x}}G^{\bf i^{\p}}J^{\bf j^{\p}}
v_{\bf k^{\p},i^{\p},j^{\p}}\\
&=(\alpha+r+k-x\beta)L^{\k^{*}}[G_{k+r-x},L^{{\k^{\p}-\k^{*}}-\es_{x}}]G^{\i}J^{\j}v_{\bf k^{\p},i^{\p},j^{\p}}\\
&+(\alpha+r+k-x\beta)L^{{\k^{\p}}-\es_{x}}G^{\i^{\p}}[G_{k+r-x},J^{\j^{\p}}]v_{\bf k^{\p},i^{\p},j^{\p}}.
\end{eqnarray*}
If it is a zero vector, there is nothing to prove.
Suppose it is nonzero.
Take $r^{'}=r-x>0$ and analyse similarly as above,
we get
${\bf w}({\k}_{x}+{\i}_{x}+{\j}_{x})\leqslant{\bf w(\k^{\p}+\i^{\p}+\j^{\p})}-r$.
Note that ${\bf w(\k^{\p}+\i^{\p}+\j^{\p})}<{\bf w(\k+\i+\j)}$,
or ${\bf w(\k^{\p}+\i^{\p}+\j^{\p})}={\bf w(\k+\i+\j)}$ and ${\bf w}({\k}_{x})<{\bf w}(\k)$.
Therefore, we have $({\k}_{x},{\i}_{x},{\j}_{x})\prec({\bf k,i,j}-\es_{r})$ in all cases.
\end{itemize}

 Combining all the arguments, we see that $\d(G_{k+r}v) =({\bf k,i,j}-\es_{r})$.

\end{proof}

\bl{lemred1.2}
Let $d\in\Z _+$ and $V$ be an $\tl_d$-module
 satisfying the conditions in Theorem \ref{idmod1}.
 For any $v\in {\rm Ind}_{\tl_{d}}^{\tl}(V)\backslash V$,
 we write $v$ in the form of (\ref{eq:idelt}).
 Suppose that ${\rm deg}(v)= ({\bf k,i,0})$
with $\i\ne 0$.
 Let $q={\rm min}\{s:i_s \ne 0\}>0$,
then ${\rm deg}(J_{k+d+q}v) =({\bf k,i}-\es_{q},{\bf0})$.
\el
\begin{proof}
 It suffices to consider those $v_{\k^{\p},\i^{\p},\j^{\p}}$ with
 $J_{k+d+q}L^{\k^{\p}}G^{\i^{\p}}J^{\j^{\p}}v_{\k^{\p},\i^{\p},\j^{\p}}\ne 0$.
Note that $J_{k+d+q}v_{\k^{\p},\i^{\p},\j^{\p}}=0$
 for any such ${(\k^{\p},\i^{\p},\j^{\p})} \in \mbox{supp}(v)$.
We have
\begin{eqnarray}
&&{}J_{k+d+q}L^{\k^{\p}}G^{\i^{\p}}J^{\j^{\p}}v_{\k^{\p},\i^{\p},\j^{\p}}=
             [J_{k+d+q},L^{\k^{\p}}]G^{\i^{\p}}J^{\j^{\p}}v_{\k^{\p},\i^{\p},\j^{\p}} \nonumber\\
&&{}\;\;\;\;+L^{\k^{\p}}[J_{k+d+q},G^{\i^{\p}}]J^{\j^{\p}}v_{\k^{\p},\i^{\p},\j^{\p}}
            +L^{\k^{\p}}G^{\i^{\p}}[J_{k+d+q},J^{\j^{\p}}]v_{\k^{\p},\i^{\p},\j^{\p}}.
\end{eqnarray}

\vspace{3mm}

{\bf Claim \ref{lemred1.2}.1:} If $\j^{\p}\neq {\bf 0}$, then
either
$L^{\k^{\p}}G^{\i^{\p}}[J_{k+d+q},J^{\j^{\p}}]v_{\k^{\p},\i^{\p},\j^{\p}}=0$
or we have
$\d(L^{\k^{\p}}G^{\i^{\p}}[J_{k+d+q},J^{\j^{\p}}]v_{\k^{\p},\i^{\p},\j^{\p}})
  \prec({\k},{\i}-\es_{q},{\bf 0})$.

{\em Proof of Claim \ref{lemred1.2}.1:}
If the central element $C_{3}$ acts as zero or $j^{\p}_{k+d+q}=0$, then we get
$L^{\k^{\p}}G^{\i^{\p}}[J_{k+d+q},J^{\j^{\p}}]v_{\k^{\p},\i^{\p},\j^{\p}}=0$.

Suppose $C_{3}$ acts as a nonzero scalar $\ell_{3}$ and $j^{\p}_{k+d+q}\neq 0$.
Since
\begin{eqnarray}
L^{\k^{\p}}G^{\i^{\p}}[J_{k+d+q},J^{\j^{\p}}]v_{\k^{\p},\i^{\p},\j^{\p}}
=(k+d+q)\ell_{3}L^{\k^{\p}}G^{\i^{\p}}J^{{\j^{\p}}-\es_{k+d+q}}v_{\k^{\p},\i^{\p},\j^{\p}}.
\end{eqnarray}
We have $\d(L^{\k^{\p}}G^{\i^{\p}}[J_{k+d+q},J^{\j^{\p}}]v_{\k^{\p},\i^{\p},\j^{\p}})
=({\k^{\p}},{\i^{\p}},{\j^{\p}}-\es_{k+d+q})$.
Then
\begin{eqnarray*}
{\bf w}({\k^{\p}}+{\i^{\p}}+{\j^{\p}}-\es_{k+d+q})
={\bf w}({\k^{\p}}+{\i^{\p}}+{\j^{\p}})-(k+d+q)
<{{\bf w}(\k+\i)}-q,
\end{eqnarray*}
which implies that
$\d(L^{\k^{\p}}G^{\i^{\p}}[J_{k+d+q},J^{\j^{\p}}]v_{\k^{\p},\i^{\p},\j^{\p}})
  \prec({\k},{\i}-\es_{q},{\bf 0})$.

\vspace{3mm}

{\bf Claim \ref{lemred1.2}.2:} If $\i^{\p}\neq {\bf 0}$, then
$\d(L^{\k^{\p}}[J_{k+d+q},G^{\i^{\p}}]J^{\j^{\p}}v_{\k^{\p},\i^{\p},\j^{\p}})
  \preceq({\k},{\i}-\es_{q},{\bf 0})$,
  where the equality holds if and only if $\k^{\p}=\k,\i^{\p}=\i, \j^{\p}={\bf 0}$.

{\em Proof of Claim \ref{lemred1.2}.2:} Suppose $J_{k+d+q}$ makes commutators with $G_{-d-y}$.
Denote
\begin{eqnarray*}
\d(L^{\k^{\p}}G^{\i^{*}}[J_{k+d+q},G_{-d-y}]G^{{\i^{\p}-\i^{*}}-{\es_{y}}}J^{\j^{\p}}v_{\k^{\p},\i^{\p},\j^{\p}})
=({\k}_{y},{\i}_{y},{\j}_{y}).
\end{eqnarray*}

\begin{itemize}
\item[(i)] For $y> k+d+q$, we have
$({\k}_{y},{\i}_{y},{\j}_{y})=({\k^{\p}},{\i^{\p}}-\es_{y}+\es_{y-q-k-d},{\j^{\p}})$.
Then
\begin{eqnarray*}
{\bf w}({\k}_{y}+{\i}_{y}+{\j}_{y})
={\bf w(\k^{\p}+\i^{\p}+\j^{\p})}-k-d-q
<{\bf w(\k+\i)}-q.
\end{eqnarray*}
\item[(ii)] For $q<y\leqslant k+d+q$,
we have
$({\k}_{y},{\i}_{y},{\j}_{y})\preceq({\k^{\p}},{\i^{\p}}-\es_{y},{\j^{\p}})$.
Hence
\begin{eqnarray*}
{\bf w}({\k}_{y}+{\i}_{y}+{\j}_{y})
\leqslant{\bf w(\k^{\p}+\i^{\p}+\j^{\p})}-y
<{\bf w(\k+\i)}-q.
\end{eqnarray*}
\item[(iii)]
 For $y=q$, we have
$({\k}_{y},{\i}_{y},{\j}_{y})=({\k^{\p}},{\i^{\p}}-\es_{q},{\j^{\p}}).$
If ${\bf w(\k^{\p}+\i^{\p}+j^{\p})}<{\bf w(\k+\i)}$, then
${\bf w}({\k}_{y}+{\i}_{y}+{\j}_{y})
={\bf w(\k^{\p}+\i^{\p}+\j^{\p})}-q<{\bf w}({\bf k+i}-\es_{q}+{\bf 0})$,
so we get $({\k}_{y},{\i}_{y},{\j}_{y})\prec({\bf k,i}-\es_{q},{\bf 0})$.
If ${\bf w(\k^{\p}+\i^{\p}+j^{\p})}={\bf w(\k+\i)}$
and ${\bf w}(\k^{\p})<{\bf w}(\k)$, then
\begin{eqnarray*}
{\bf w}({\k}_{y}+{\i}_{y}+{\j}_{y})
={\bf w(\k^{\p}+\i^{\p}+\j^{\p})}-q={\bf w}({\bf k+i}-\es_{q}+{\bf 0}),
\end{eqnarray*}
but ${\bf w}({\k}_{y})={\bf w(\k^{\p})}<{\bf w}(\k)$.
So we also get
$({\k}_{y},{\i}_{y},{\j}_{y})\prec({\bf k,i}-\es_{q},{\bf 0})$.
If ${\bf w(\k^{\p}+\i^{\p}+j^{\p})}={\bf w(\k+\i+\j)}$
and ${\bf w}(\k^{\p})={\bf w}(\k)$,
then we must have $\j^{\p}=0$.
Therefore,
$({\k}_{y},{\i}_{y},{\j}_{y})\preceq({\bf k,i}-\es_{q},{\bf 0})$,
and the equality holds if and only if ${\bf k^{\p}=k,i^{\p}=i,j^{\p}=0}$.
\item[(iv)] For $1\leqslant y<q$,
we have
${\bf w}({\k}_{y}+{\i}_{y}+{\j}_{y})\leqslant{\bf w(\k^{\p}+\i^{\p}+\j^{\p})}-q.$
If ${\bf w(\k^{\p}+\i^{\p}+\j^{\p})}<{\bf w(\k+\i)}$, then
$({\k}_{y},{\i}_{y},{\j}_{y})\prec({\bf k,i}-\es_{q},{\bf 0})$.
If ${\bf w(\k^{\p}+\i^{\p}+\j^{\p})}={\bf w(\k+\i)}$,
then ${\bf w}({\k^{\p}})<{\bf w}(\k)$.
So ${\bf w}({\k}_{y})={\bf w}({\k^{\p}})<{\bf w}(\k)$.
Hence, $({\k}_{y},{\i}_{y},{\j}_{y})\prec({\bf k,i}-\es_{q},{\bf 0})$.
\end{itemize}

\vspace{3mm}

{\bf Claim \ref{lemred1.2}.3:} If $\k^{\p}\neq {\bf 0}$, then
$\d([J_{k+d+q},L^{\k^{\p}}]G^{\i^{\p}}J^{\j^{\p}}v_{\k^{\p},\i^{\p},\j^{\p}})
  \prec({\k},{\i}-\es_{q},{\bf 0})$.

{\em Proof of Claim \ref{lemred1.2}.3:} Suppose $J_{k+d+q}$ makes commutators with $L_{-x}$.
Denote
\begin{eqnarray*}
\d(L^{\k^{*}}[J_{k+d+q},L_{-x}]L^{{\k^{\p}-\k^{*}}-{\es_{x}}}G^{\i^{\p}}J^{\j^{\p}}v_{\k^{\p},\i^{\p},\j^{\p}})
=({\k}_{x},{\i}_{x},{\j}_{x}).
\end{eqnarray*}

 \begin{itemize}
\item[(i)] For $x> k+d+q$, we have
$({\k}_{x},{\i}_{x},{\j}_{x})=({\k^{\p}}-\es_{x},{\i^{\p}},{\j^{\p}}+\es_{x-q-k-d})$.
Then
\begin{eqnarray*}
{\bf w}({\k}_{x}+{\i}_{x}+{\j}_{x})
={\bf w(\k^{\p}+\i^{\p}+\j^{\p})}-k-d-q
<{\bf w(\k+\i)}-q.
\end{eqnarray*}
\item[(ii)] For $q\leqslant x\leqslant k+d+q$,
we get that
$({\k}_{x},{\i}_{x},{\j}_{x})\preceq({\k^{\p}}-\es_{x},{\i^{\p}},\j^{\p})$.
If ${\bf w(\k^{\p}+\i^{\p}+j^{\p})}<{\bf w(\k+\i)}$, then
${\bf w}({\k}_{x}+{\i}_{x}+{\j}_{x})
\leqslant{\bf w(\k^{\p}+\i^{\p}+\j^{\p})}-x<{\bf w}({\bf k+i}-\es_{q}+{\bf 0}).$
If ${\bf w(\k^{\p}+\i^{\p}+j^{\p})}={\bf w(\k+\i)}$, then
${\bf w}({\k}_{x}+{\i}_{x}+{\j}_{x})
\leqslant{\bf w(\k^{\p}+\i^{\p}+\j^{\p})}-q={\bf w}({\bf k+i}-\es_{q}+{\bf 0}),$
but ${\bf w} ({\k}_{x})={\bf w(\k^{\p})}-x<{\bf w}(\k)$.
So we always have
$({\k}_{x},{\i}_{x},{\j}_{x})\prec({\bf k,i}-\es_{q},{\bf 0})$.
\item[(iii)] For $1\leqslant x<q$,
take $q^{'}=q-x>0$ and analyse as above,
we get
$
{\bf w}({\k}_{x}+{\i}_{x}+{\j}_{x})\leqslant{\bf w(\k^{\p}+\i^{\p}+\j^{\p})}-q.
$
Since ${\bf w(\k^{\p}+\i^{\p}+\j^{\p})}<{\bf w(\k+\i+\j)}$,
or ${\bf w(\k^{\p}+\i^{\p}+\j^{\p})}={\bf w(\k+\i+\j)}$ and ${\bf w}({\k}_{x})<{\bf w}(\k)$.
We have $({\k}_{x},{\i}_{x},{\j}_{x})\prec({\bf k,i}-\es_{q},{\bf 0})$ in all cases.
\end{itemize}

All the above arguments tell that $\d(J_{k+d+q}v) =({\bf k,i}-\es_{q},{\bf 0})$.

\end{proof}

\bl{lemred1.3}
Let $d\in\Z _+$ and $V$ be an $\tl_d$-module
 satisfying the conditions in Theorem \ref{idmod1}.
 For any $v\in {\rm Ind}_{\tl_{d}}^{\tl}(V)\backslash V$,
 we write $v$ in the form of (\ref{eq:idelt}).
 Suppose ${\rm deg}(v)= ({\bf k,0,0})$
 with $\k\ne 0$.
  Let $p={\rm min}\{s:k_s \ne 0\}>0$,
then ${\rm deg}(L_{l+p}v) =({\bf k}-\es_{p},{\bf 0,0})$.
\el
 \begin{proof}
 It suffices to consider those $v_{\k^{\p},\i^{\p},\j^{\p}}$ with
 $L_{l+p}L^{\k^{\p}}G^{\i^{\p}}J^{\j^{\p}}v_{\k^{\p},\i^{\p},\j^{\p}}\ne 0$.
 Noticing that $L_{l+p}v_{\k^{\p},\i^{\p},\j^{\p}}=0$
  for any such ${(\k^{\p},\i^{\p},\j^{\p}) }\in \mbox{supp}(v)$.
  We have
\begin{eqnarray}
&&{}L_{l+p}L^{\k^{\p}}G^{\i^{\p}}J^{\j^{\p}}v_{\k^{\p},\i^{\p},\j^{\p}}=
             [L_{l+p},L^{\k^{\p}}]G^{\i^{\p}}J^{\j^{\p}}v_{\k^{\p},\i^{\p},\j^{\p}} \nonumber\\
&&{}\;\;\;\;+L^{\k^{\p}}[L_{l+p},G^{\i^{\p}}]J^{\j^{\p}}v_{\k^{\p},\i^{\p},\j^{\p}}
            +L^{\k^{\p}}G^{\i^{\p}}[L_{l+p},J^{\j^{\p}}]v_{\k^{\p},\i^{\p},\j^{\p}}.
\end{eqnarray}

\vspace{3mm}

{\bf Claim \ref{lemred1.3}.1:} If $\j^{\p}\neq {\bf 0}$, then
$\d(L^{\k^{\p}}G^{\i^{\p}}[L_{l+p},J^{\j^{\p}}]v_{\k^{\p},\i^{\p},\j^{\p}})
  \prec({\k}-\es_{p},{\bf 0},{\bf 0})$.

{\em Proof of Claim \ref{lemred1.3}.1:}
Suppose $L_{l+p}$ makes commutator with $J_{-z}$.
Denote
\begin{eqnarray*}
\d(L^{\k^{\p}}G^{\i^{\p}}J^{\j^{*}}[L_{l+p},J_{-z}]J^{{\j^{\p}-\j^{*}}-{\es_{z}}}v_{\k^{\p},\i^{\p},\j^{\p}})
=({\k}_{z},{\i}_{z},{\j}_{z}).
\end{eqnarray*}

\begin{itemize}
\item[(i)] For $z> l+p$, we have
$({\k}_{z},{\i}_{z},{\j}_{z})=({\k^{\p}},{\i^{\p}},{\j^{\p}}-\es_{z}+\es_{z-l-p})$.
Then
\begin{eqnarray*}
{\bf w}({\k}_{z}+{\i}_{z}+{\j}_{z})
={\bf w(\k^{\p}+\i^{\p}+\j^{\p})}-l-p
<{\bf w(\k)}-p.
\end{eqnarray*}
\item[(ii)] For $p+l-k-d\leqslant z\leqslant l+p$,
we have
$({\k}_{z},{\i}_{z},{\j}_{z})\preceq({\k^{\p}},{\i^{\p}},{\j^{\p}}-\es_{z})$.
Then (note that $l>k+d$)
\begin{eqnarray*}
{\bf w}({\k}_{z}+{\i}_{z}+{\j}_{z})
\leqslant{\bf w(\k^{\p}+\i^{\p}+\j^{\p})}-z
<{\bf w(\k+\i)}-p.
\end{eqnarray*}
\item[(iii)] For $1\leqslant z<p+l-k-d$.
Note that $p+l-z>0$ and
\begin{eqnarray*}
L^{\k^{\p}}G^{\i^{\p}}J^{\j^{*}}[L_{l+p},J_{-z}]J^{{\j^{\p}-\j^{*}}-{\es_{z}}}v_{\k^{\p},\i^{\p},\j^{\p}}
=zL^{\k^{\p}}G^{\i^{\p}}J^{\j^{*}}[J_{l+p-z},J^{{\j^{\p}-\j^{*}}-{\es_{z}}}]v_{\k^{\p},\i^{\p},\j^{\p}}.
\end{eqnarray*}
If it is a zero vector, there is nothing to prove.
Suppose it is nonzero,
then $({\k}_{z},{\i}_{z},{\j}_{z})=({\k^{\p}},{\i^{\p}},{\j^{\p}}-\es_{z}-\es_{l+p-z})$.
Since
${\bf w}({\k}_{z}+{\i}_{z}+{\j}_{z})\leqslant
{\bf w(\k^{\p}+\i^{\p}+\j^{\p})}-p-l
<{\bf w(\k)}-p.$
We get $({\k}_{z},{\i}_{z},{\j}_{z})\prec({\bf k}-\es_{p},{\bf 0},{\bf 0})$.
\end{itemize}

\vspace{3mm}

{\bf Claim \ref{lemred1.3}.2:} If $\i^{\p}\neq {\bf 0}$, then
$\d(L^{\k^{\p}}[L_{l+p},G^{\i^{\p}}]J^{\j^{\p}}v_{\k^{\p},\i^{\p},\j^{\p}})
  \prec({\k}-\es_{p},{\bf 0},{\bf 0})$.

{\em Proof of Claim \ref{lemred1.3}.2:}
Suppose $L_{l+p}$ makes commutator with $G_{-d-y}$.
Denote
\begin{eqnarray*}
\d(L^{\k^{\p}}G^{\i^{*}}[L_{l+p},G_{-d-y}]G^{{\i^{\p}-\i^{*}}-{\es_{y}}}J^{\j^{\p}}v_{\k^{\p},\i^{\p},\j^{\p}})
=({\k}_{y},{\i}_{y},{\j}_{y}).
\end{eqnarray*}

\begin{itemize}
\item[(i)] For $y> l+p$, we have
$({\k}_{y},{\i}_{y},{\j}_{y})=({\k^{\p}},{\i^{\p}}-\es_{y}+\es_{y-l-p},{\j^{\p}})$.
Then
${\bf w}({\k}_{y}+{\i}_{y}+{\j}_{y})
={\bf w(\k^{\p}+\i^{\p}+\j^{\p})}-l-p
<{\bf w(\k+\i)}-p.$
\item[(ii)] For $p+l-k-d\leqslant y\leqslant l+p$,
we have
$({\k}_{y},{\i}_{y},{\j}_{y})\preceq({\k^{\p}},{\i^{\p}}-\es_{y},{\j^{\p}})$.
Note that $l>k+d$. Therefore,
${\bf w}({\k}_{y}+{\i}_{y}+{\j}_{y})
\leqslant{\bf w(\k^{\p}+\i^{\p}+\j^{\p})}-y
<{\bf w(\k+\i)}-p.$
\item[(iii)] For $1\leqslant y<p+l-k-d$,
we have
${\bf w}({\k}_{y}+{\i}_{y}+{\j}_{y})\leqslant
{\bf w(\k^{\p}+\i^{\p}+\j^{\p})}-p-(l-k-d)
<{\bf w(\k)}-p.$
\end{itemize}

Hence, $({\k}_{y},{\i}_{y},{\j}_{y})\prec({\bf k}-\es_{p},{\bf 0},{\bf 0})$.

\vspace{3mm}

{\bf Claim \ref{lemred1.3}.3:} If $\k^{\p}\neq {\bf 0}$, then
$\d([L_{l+p},L^{\k^{\p}}]G^{\i^{\p}}J^{\j^{\p}}v_{\k^{\p},\i^{\p},\j^{\p}})
  \preceq({\k}-\es_{p},{\bf 0},{\bf 0})$.
And the equality holds if and only if $\k^{\p}=\k,\i^{\p}=\j^{\p}={\bf 0}$.

{\em Proof of Claim \ref{lemred1.3}.3:}
Suppose that $L_{l+p}$ makes commutator with $L_{-x}$.
Denote the degree of the vector by
\begin{eqnarray*}
\d(L^{\k^{*}}[L_{l+p},L_{-x}]L^{{\k^{\p}-\k^{*}}-\es_{x}}G^{\i^{\p}}J^{\j^{\p}}v_{\k^{\p},\i^{\p},\j^{\p}})
=({\k}_{x},{\i}_{x},{\j}_{x}).
\end{eqnarray*}

 \begin{itemize}
\item[(i)] For $x> l+p$, we have
$({\k}_{x},{\i}_{x},{\j}_{x})=({\k^{\p}}-\es_{x}+\es_{x-l-p},{\i^{\p}},{\j^{\p}})$.
Then (note that $l>0$)
\begin{eqnarray*}
{\bf w}({\k}_{x}+{\i}_{x}+{\j}_{x})
={\bf w(\k^{\p}+\i^{\p}+\j^{\p})}-l-p
<{\bf w(\k)}-p.
\end{eqnarray*}
\item[(ii)] For $p< x\leqslant l+p$,
we have
$({\k}_{x},{\i}_{x},{\j}_{x})\preceq({\k^{\p}}-\es_{x},{\i^{\p}},\j^{\p})$.
And then
${\bf w}({\k}_{x}+{\i}_{x}+{\j}_{x})
\leqslant{\bf w(\k^{\p}+\i^{\p}+\j^{\p})}-x<{\bf w}({\k})-p.$
\item[(iii)] For $x=p$, we have
$
({\k}_{x},{\i}_{x},{\j}_{x})=({\k^{\p}}-\es_{p},{\i^{\p}},{\j^{\p}}).
$
If ${\bf w}({\k}^{\p}+{\i}^{\p}+{\j}^{\p})<{\bf w}(\k)$, then
${\bf w}({\k}_{x}+{\i}_{x}+{\j}_{x})
={\bf w}({\k}^{\p}+{\i}^{\p}+{\j}^{\p})-p<{\bf w(\k)}-p,$
so that $({\k}_{x},{\i}_{x},{\j}_{x})\prec({\k}-\es_{p},{\bf 0},{\bf 0})$.
If ${\bf w}({\k}^{\p}+{\i}^{\p}+{\j}^{\p})={\bf w}(\k)$
and ${\bf w}(\k^{\p})<{\bf w}(\k)$,
then
$
{\bf w}({\k}_{x}+{\i}_{x}+{\j}_{x})
={\bf w}({\k}^{\p}+{\i}^{\p}+{\j}^{\p})-p={\bf w(\k)}-p,
$
but ${\bf w}({\k}_{x})={\bf w(\k^{\p})}-p<{\bf w(\k)}-p$,
so that $({\k}_{x},{\i}_{x},{\j}_{x})\prec({\k}-\es_{p},{\bf 0},{\bf 0})$.
If ${\bf w}({\k}^{\p}+{\i}^{\p}+{\j}^{\p})={\bf w}(\k)$
and ${\bf w}(\k^{\p})={\bf w}(\k)$, then
we must have $\i^{\p}=\j^{\p}={\bf 0}$
and $\k^{\p}\preceq\k$.
Hence,
$({\k}_{x},{\i}_{x},{\j}_{x})\preceq({\k}-\es_{p},{\bf 0},{\bf 0}),$
where the equality holds if and only if $\k^{\p}=\k,\i^{\p}=\j^{\p}={\bf 0}$.
\item[(iv)] For $1\leqslant x<p$,
we have
$
{\bf w}({\k}_{x}+{\i}_{x}+{\j}_{x})\leqslant{\bf w(\k^{\p}+\i^{\p}+\j^{\p})}-p.
$
Note that ${\bf w(\k^{\p}+\i^{\p}+\j^{\p})}<{\bf w(\k+\i+\j)}$,
or ${\bf w(\k^{\p}+\i^{\p}+\j^{\p})}={\bf w(\k+\i+\j)}$ and ${\bf w}({\k}_{x})<{\bf w}(\k)$.
Therefore, we have $({\k}_{x},{\i}_{x},{\j}_{x})\prec({\k}-\es_{p},{\bf 0},{\bf 0})$.
\end{itemize}

 In summary, we have $\d(L_{k+p}v) =({\k}-\es_{p},{\bf 0},{\bf 0})$.
 \end{proof}

\vspace{3mm}

{\bf 3.2. The case of $F=0$, $C_{3}\neq 0$.}

\vspace{3mm}

In this subsection, we give a construction of simple restricted $\tl$-modules
for $F=0$, where $C_{3}$ acts as a nonzero scalar (we will write as $C_{3}\neq 0$ for brevity).
Here we deal with the cases of $\alpha,\beta$ such that $\alpha+n+m\beta\neq 0$
for $n,m\in\Z\backslash\{0\}$.

For any $d_1,d_2\in \Z_+$,
 let $\ud=(d_1,d_2)$.
 Consider the subalgebra
\begin{equation}
\tl_{\ud}=\sum_{i\in\Z_+}(\C G_{-d_1+i}\oplus\C J_{-d_2+i}\oplus\C L_i)\oplus \sum_{i=1}^{3}\C C_i.
\end{equation}
Let $V$ be a simple $\tl_{\ud}$-module with $C_i$ acts as
a scalar for $i=1,2,3$.
Form the induced $\tl$-module
\begin{equation}\nonumber
\id_{\tl_{\ud}}^{\tl}(V) =U(\tl)\otimes _{U(\tl_{\ud})}V.
\end{equation}

We show that under certain conditions, the induced module is a simple $\tl$-module.
\bt{idmod2}
Let $F=0$ and $C_3\ne0$.
    Assume that $\ud=(d_1,d_2)\in\Z_+\times\Z_{+}$
    and $V$ is a simple $\tl_{\ud}$-module.
    If there exists $k\in\Z_+$ such that
\vspace{3mm}\\
(a) $G_{k}$ acts injectively on $V$,\\
(b) $G_{k+i}V=J_{d_{2}+i}V=L_{k+d_{1}+i}V=0$ for all $i\in\N$,
\vspace{2mm}\\
then ${\rm Ind}_{\tl_{\ud}}^{\tl}(V)$ is a simple $\tl$-module.
\et

Before giving the proof, we introduce some notions and notations that will be used.
 For ${\a,\b,\c}\in\M$, denote
\begin{equation}\nonumber
G^{\a}J^{\b}L^{\c}
=\cdots G_{-d_1-2}^{a_2}G_{-d_1-1}^{a_1}\cdots J_{-d_2-2}^{b_2}J_{-d_2-1}^{b_1}
\cdots L_{-2}^{c_2}L_{-1}^{c_1}\in U(\tl).
\end{equation}
By the PBW Theorem, every element of $\id_{\tl_{\ud}}^{\tl}(V)$ can be uniquely written in the form
\begin{equation}\label{eq:idelt2}
  \sum_{{\a,\b,\c}\in\M}G^{\a}J^{\b}L^{\c}v_{\a,\b,\c},
\end{equation}
where all $v_{\a,\b,\c}\in V$ and only finitely many of them are nonzero.

\bd{pto2}
{\em Define a {\em principal total order} $({\bf o'})$ on $\M\times\M\times\M$,
 still denoted by $\succ$, as follows:
 for any $ {\bf a,b,c,a^{\p},b^{\p},c^{\p}} \in \M$,
 we say that
 $({\bf a,b,c}) \succ (\a^{\p},\b^{\p},\c^{\p})$
 if one of the following conditions hold.
  \flushleft(i)   $({\c},{\bf w(c)})\succ({\c^{\p}},{\bf w(c^{\p})})$;\\
 (ii)  $\c=\c^{\p}, \b\succ \b^{\p}$;\\
 (iii) $\c=\c^{\p},\b=\b^{\p},\a>\a^{\p}$.
 }
\ed

For any $v\in\id_{\tl_{\ud}}^{\tl}(V)$ written in the form of (\ref{eq:idelt2}),
we denote by $\mbox{supp}(v)$ the set of all $({\bf a,b,c})\in \M\times\M\times\M$
such that $v_{\bf a,b,c}\ne 0$.
 For a nonzero $v\in\id_{\tl_{\ud}}^{\tl}(V)$,
   let $\d(v)$ denote the maximal
    (with respect to the principal total order $({\bf o'})$)
     element in $\mbox{supp}(v)$,
      called the {\em degree} of $v$.

\bl{lemred2}
Let $\ud=(d_1,d_2)\in\Z_+\times\Z_{+}$ and
$V$ be an $\tl_{\ud}$-module
 satisfying the conditions in Theorem \ref{idmod2}.
  For any $v\in{\rm Ind}_{\tl_{\ud}}^{\tl}(V)\backslash V$
  written in the form of (\ref{eq:idelt2}),
  denote ${\rm deg}(v)= ({\bf a,b,c})$.
  \vspace{3mm}\\
 {\rm (1)} If ${\bf c\ne 0}$, $r={\rm min}\{s:c_s \ne 0\}>0$,
 then ${\rm deg}(G_{k+r}v) =({\bf a,b},{\c}-\es_{r})$.\\
{\rm (2)} If ${\bf c=0,b\ne 0}$,
$q={\rm min}\{s:b_s \ne 0\}>0$,
then ${\rm deg}(J_{d_2+q}v) =({\a},{\b}-\es_{q},{\bf0})$.\\
{\rm (3)} If ${\bf c=b=0,a\ne 0}$,
$p={\rm max}\{s:a_s \ne 0\}>0$,
then ${\rm deg}(L_{k+d_1+p}v) =({\a}-\es_{p},{\bf0,0})$.
\el
\begin{proof}
(1) It suffices to consider those $v_{\bf a^{\p},b^{\p},c^{\p}}$
with $G_{k+r}G^{\bf a^{\p}}J^{\bf b^{\p}}L^{\bf c^{\p}}v_{\bf a^{\p},b^{\p},c^{\p}}\ne 0$.
Noticing that $G_{k+r}v_{\bf a^{\p},b^{\p},c^{\p}}=0$
for any such $({\bf a^{\p},b^{\p},c^{\p}}) \in \mbox{supp}(v)$.
We have
\begin{eqnarray}
G_{k+r}G^{\bf a^{\p}}J^{\bf b^{\p}}L^{\bf c^{\p}}v_{\bf a^{\p},b^{\p},c^{\p}}
=G^{\bf a^{\p}}J^{\bf b^{\p}}[G_{k+r},L^{\bf c^{\p}}]v_{\bf a^{\p},b^{\p},c^{\p}}.
\end{eqnarray}

Suppose $G_{k+r}$ makes commutator with $L_{-x}$.
Denote
\begin{eqnarray*}
\d(G^{\bf a^{\p}}J^{\bf b^{\p}}L^{\bf c^{*}}[G_{k+r},L_{-x}]L^{{\c^{\p}-\c^{*}}-\es_{x}}v_{\bf a^{\p},b^{\p},c^{\p}})
=({\a}_{x},{\b}_{x},{\c}_{x}).
\end{eqnarray*}
Note that by our assumption of $\alpha$ and $\beta$, $[G_{k+r},L_{-x}]\neq 0$.

\begin{itemize}
\item[(i)] For $x>k+d+r$, we have
$({\a}_{x},{\b}_{x},{\c}_{x})=({\a^{\p}}+\es_{x-k-d-r},{\b^{\p}},{\c^{\p}}-\es_{x})$.
Since ${\bf w({\c}}_{x})={\bf w}({\c^{\p}}-\es_{x})={\bf w(\c^{\p})}-x<{\bf w(\c)}-r$,
we get $({\a}_{x},{\b}_{x},{\c}_{x})\prec({\bf a,b},{\c}-\es_{r})$.
\item[(ii)] For $r<x\leqslant k+d+r$, we have
$({\a}_{x},{\b}_{x},{\c}_{x})\preceq({\a^{\p}},{\b^{\p}},{\c^{\p}}-\es_{x})$.
It follows from
 ${\bf w({\c}}_{x})\leqslant{\bf w}({\c^{\p}}-\es_{x})={\bf w(\c^{\p})}-x<{\bf w(\c)}-r$
that $({\a}_{x},{\b}_{x},{\c}_{x})\prec({\bf a,b},{\c}-\es_{r})$.
\item[(iii)] For $x=r$, we see that
$({\a}_{x},{\b}_{x},{\c}_{x})=({\a^{\p}},{\b^{\p}},{\c^{\p}}-\es_{r})$.
If ${\bf w}(\c^{\p})<{\bf w}(\c)$, then
${\bf w}({\c}_{x})={\bf w(\c^{\p})}-r<{\bf w(\c)}-r$
tells that $({\a}_{x},{\b}_{x},{\c}_{x})\prec({\bf a,b},{\c}-\es_{r})$.
If ${\bf w}(\c^{\p})={\bf w}(\c)$, we have $\c^{\p}\preceq\c$.
Then
$({\a}_{x},{\b}_{x},{\c}_{x})=({\a^{\p}},{\b^{\p}},{\c^{\p}}-\es_{r})\preceq({\bf a,b},{\c}-\es_{r}),$
where the equality holds if and only if $\a^{\p}=\a,\b^{\p}=\b,\c^{\p}=\c$.
\item[(iv)] For $1\leq x<r$, we have ${\bf w(\c^{\p})< w(\c)}$.
Then
${\bf w}({\c}_{x})\leqslant{\bf w(\c^{\p})}-r<{\bf w(\c)}-r.$
Hence, $({\a}_{x},{\b}_{x},{\c}_{x})\prec({\bf a,b},{\c}-\es_{r})$.
\end{itemize}

Therefore, we have ${\rm deg}(G_{k+r}v) =({\bf a,b},{\c}-\es_{r})$.

\vspace{4mm}

(2)
If $\c=0,\b\neq0$, then
\begin{eqnarray}
v=\sum\limits_{(\a^{\p},\b^{\p},0)\preceq(\a,\b,0)}G^{\a^{\p}}J^{\b^{\p}}v_{\a^{\p},\b^{\p},0}.
\end{eqnarray}
Since $J_{d_{2}+q}v_{\a^{\p},\b^{\p},0}=0$, we have (note that $C_{3}\neq0$)
\begin{eqnarray*}
J_{d_{2}+q}G^{\a^{\p}}J^{\b^{\p}}v_{\a^{\p},\b^{\p},0}
=G^{\a^{\p}}[J_{d_{2}+q},J^{\b^{\p}}]v_{\a^{\p},\b^{\p},0}.
\end{eqnarray*}

Let $q_{1}={\rm min}\{s: b^{\p}_s \ne 0\}$.
Then $q_{1}\geq q$.
If $q_{1}>q$, then $J_{d_{2}+q}G^{\a^{\p}}J^{\b^{\p}}v_{\a^{\p},\b^{\p},0}=0$.
For $q_{1}=q$, we have
\begin{eqnarray*}
J_{d_{2}+q}G^{\a^{\p}}J^{\b^{\p}}v_{\a^{\p},\b^{\p},0}
=G^{\a^{\p}}J^{{\b^{\p}}-\es_{q}}[J_{d_2+q},J_{-d_{2}-q}]v_{\a^{\p},\b^{\p},0}
=(d_{2}+q)C_{3}G^{\a^{\p}}J^{{\b^{\p}}-\es_{q}}v_{\a^{\p},\b^{\p},0}.
\end{eqnarray*}
Hence,
$\d(J_{d_{2}+q}G^{\a^{\p}}J^{\b^{\p}}v_{\a^{\p},\b^{\p},0})
=(\a^{\p},{\b^{\p}}-\es_{q},{\bf 0})\preceq(\a,{\b}-\es_{q},{\bf 0})$,
where the equality holds if and only if $\a^{\p}=\a,\b^{\p}=\b$.

\vspace{4mm}

(3)
Suppose $\b=\c=0, \a\neq0$, then
\begin{eqnarray}
v=\sum\limits_{\a^{\p}\leqslant\a}G^{\a^{\p}}v_{\a^{\p},0,0}.
\end{eqnarray}
Note that $L_{k+d_{1}+p}V=0$ and
\begin{eqnarray*}
[L_{k+d_{1}+p},G_{-d_{1}-p}]V\neq0,\;\;
[L_{k+d_{1}+p},G_{-d_{1}-a}]V=0\;\;\mbox{for}\;a<p.
\end{eqnarray*}
Let $p_{1}={\rm max}\{s: a^{\p}_{s}\neq 0\}$, then
$p_{1}\leqslant p$.
If $p_{1}<p$, then $L_{k+d_{1}+p}G^{\a^{\p}}v_{\a^{\p},0,0}=0$.
For $p_{1}=p$,
we have
\begin{eqnarray*}
&L_{k+d_{1}+p}G^{\a^{\p}}v_{\a^{\p},0,0}
=G^{{\a^{\p}}-\es_{p}}[L_{k+d_{1}+p},G_{-d_{1}-p}]v_{\a^{\p},0,0}\\
&\;\;=-(\alpha-d_{1}-p+(k+d_{1}+p)\beta)G^{{\a^{\p}}-\es_{p}}G_{k}v_{\a^{\p},0,0}.
\end{eqnarray*}
Since we require that $(\alpha-(d_{1}+p)+(k+d_{1}+p)\beta)\neq 0$,
we get
$\d(L_{k+d_{1}+p}G^{\a^{\p}}v_{\a^{\p},0,0})
=({\a^{\p}}-\es_{p},{\bf0,0})\preceq({\a}-\es_{p},{\bf0,0})$,
and the equality holds if and only if $\a^{\p}=\a$.
\end{proof}

Using Lemma \ref{lemred2} repeatedly,
from any $0\neq v\in \id_{\tl_{\ud}}^{\tl}(V)$
we can reach a nonzero element in
$U(\tl)v\cap V \ne 0$,
which gives the simplicity of $\id_{\tl_{\ud}}^{\tl}(V)$.

\section{Characterization of simple restricted $\tl$-modules}
\label{sec:4}
\def\theequation{4.\arabic{equation}}
	\setcounter{equation}{0}

In this section, we present a precise characterization of simple
restricted $\tl$-modules under certain conditions.
We assume that
$\alpha+n+m\beta\neq 0$ for $m\in\Z,n\in\Z\backslash\{0\}$ as needed.
For any $k,l,m\in\Z_{+}$,
denote
\begin{eqnarray}
{\tl}^{(k,l,m)}=\sum\limits_{i\in\Z_{+}}(\C G_{k+i}\oplus\C J_{l+i}\oplus\C L_{m+i}).
\end{eqnarray}
It is straightforward to check that ${\tl}^{(k,l,m)}$ is a finitely generated subalgebra of $\tl$.

First, we show several equivalent
conditions for simple restricted modules over $\tl$.
\bt{eqchar}
  Suppose that $S$ is a simple $\tl$-module.
  Then the following conditions are equivalent:
  \vspace{3mm}\\
{\rm (1)} There exists $k\in\Z_{+}$ such that $G_{i},J_i,L_i, i\geqslant k$ act on $S$ locally finitely.\\
{\rm (2)} There exists $k\in\Z_{+}$ such that $G_{i},J_i,L_i, i\geqslant k$ act on $S$ locally nilpotently.\\
{\rm (3)} There exist $k,l,m\in\Z_{+}$ such that $S$ is a locally finite ${\tl}^{(k,l,m)}$-module.\\
{\rm (4)} There exist $k,l,m\in\Z_{+}$ such that $S$ is a locally nilpotent ${\tl}^{(k,l,m)}$-module.\\
{\rm (5)} There exist $k,l,m\in\Z_{+}$ and a nonzero vector $v\in S$ such that ${\tl}^{(k,l,m)}v=0$;\\
{\rm (6)} $S$ is restricted.
\et
\begin{proof}
The conclusions
$(4)\Rightarrow(2)\Rightarrow(1)$,
$(3)\Rightarrow(1)$ and $(6)\Rightarrow(5)$ are clear.
Since ${\tl}^{(k,l,m)}$ is a finitely generated Lie algebra,
we have $(4)\Rightarrow(3)$.
Suppose $(5)$ holds.
Since $S$ is a simple $\tl$-module,
we have $S=U(\tl)v$ for the vector $v$ in (5).
Then by the PBW Theorem and the Lie brackets of $\tl$
we deduce that $S$
is a locally nilpotent ${\tl}^{(k^{\p},l^{\p},m^{\p})}$-module
for some $k^{\p},l^{\p},m^{\p}\in\Z_{+}$.
So we get
$(5)\Rightarrow(4)$.
Therefore, we only need to show that $(1)\Rightarrow(5)$
and $(5)\Rightarrow(6)$.

$(1)\Rightarrow(5)$:
Since $L_k$ acts locally finitely on $S$,
 there exists a nonzero $v\in S$ such that $L_kv=\lambda v$ for some $\lambda\in\C$.
For any $j\in\Z$ with $j>k$, set
\begin{equation}
L(j)=\sum_{n\in\Z _+}{\C L_{k}^{n}}L_jv,\;\;\;
J(j)=\sum_{n\in \Z_+}{\C L_{k}^{n}}J_jv,\;\;\;
G(j)=\sum_{n\in \Z_+}{\C L_{k}^{n}}G_{j}v,
\end{equation}
which are all finite-dimensional.
We have
\begin{equation}
\begin{aligned}
&-(j+( n-1)k) L_{j+(n+1)k}v
=[L_k,L_{j+nk}] v=( L_k-\lambda) L_{j+nk}v,
\\&-(j+nk)J_{j+(n+1)k}v=[L_k,J_{j+nk}]v
=(L_k-\lambda)J_{j+nk}v,
\\&-(\alpha+j+(n+\beta)k)G_{j+(n+1)k}v=[L_k,G_{j+nk}]v
=(L_k-\lambda) G_{j+nk}v
\end{aligned}
\end{equation}
for all $n\in\Z_+$.
Note that we assume
$\alpha+j+(n+\beta)k\neq 0$.
Therefore,
 $L_{j+nk}v\in L(j)$,
 $J_{j+nk}v\in J(j)$,
 $G_{j+nk}v\in G(j)$
  imply that
 $L_{j+(n+1)k}v\in L(j)$,
 $J_{j+(n+1)k}v\in J(j)$,
 $G_{j+(n+1) k}v\in G(j)$ respectively.
By induction on $n$, we get
$L_{j+nk}v\in L(j),J_{j+nk}v\in J(j),G_{j+nk}v\in G(j)$
for all $n\in \Z _+$.
In particular, $\sum\limits_{n\in \Z_+}{\C L_{j+nk}v}$,
$\sum\limits_{n\in \Z_+}{\C J_{j+nk}v}$
and $\sum\limits_{n\in\Z_+}{\C G_{j+nk}v}$ are finite-dimensional for all $j>k$.
Hence,
\begin{equation}
\begin{aligned}
&\sum_{i\in \Z_+}{\C L_{k+i}v}=\C L_kv+\sum_{j=k+1}^{2k}{(\sum_{n\in\Z_+}{\C L_{j+nk}v})},
\\&\sum_{i\in\Z_+}{\C J_{k+i}v}=\C J_kv+\sum_{j=k+1}^{2k}{(\sum_{n\in\Z_+}{\C J_{j+nk}v})},
\\&\sum_{i\in\Z_+}{\C {G_{k+i}}v}=\C {G_{k}}v+\sum_{j=k+1}^{2k}{(\sum_{n\in\Z_+}{\C {G_{j+nk}}v})}
\end{aligned}
\end{equation}
are all finite-dimensional.
So we can find $l\in\N$ such that
\begin{equation}\label{eq:1}
\sum_{i\in\Z_+}{\C L_{k+i}v}=\sum_{i=0}^l{\C L_{k+i}v},\;\;
\sum_{i\in\Z_+}{\C J_{k+i}v}=\sum_{i=0}^l{\C J_{k+i}v},\;\;
\sum_{i\in\Z_+}{\C G_{k+i}v}=\sum_{i=0}^l{\C G_{k+i}v}.
\end{equation}
Denote
$$V'=\sum\limits_{k_i,n_i,m_i\in\Z_+,
\\ i=0,\ldots,l}{\C G_{k}^{k_0}\cdots G_{k+l}^{k_l}}J_{k}^{n_0}\cdots J_{k+l}^{n_l}L_{k}^{m_0}\cdots L_{k+l}^{m_l}v$$
 which is nonzero and finite-dimensional by condition (1).
Using Lie brackets of $\tl$ and (\ref{eq:1})
we show that $V'$ is an $\tl^{(k,k,k)}$-module.

Let $s\in\Z_+$ be minimal such that
$(L_m+a_1L_{m+1}+\cdots +a_sL_{m+s})V'=0$
 for some $m\geqslant k$ and $a_i\in\C$. Applying $L_m$, we get
\begin{equation*}
(a_1[ L_m,L_{m+1}] +\cdots +a_s[ L_m,L_{m+s}])V'=0.
\end{equation*}
To avoid a contradiction, we must have $s=0$, that is, $L_mV'=0$.
Then, for any $i\geqslant k$, we have
\begin{equation*}
0=L_iL_mV'=[L_i,L_m] V'+L_mL_iV'=(i-m) L_{m+i}V'.
\end{equation*}
Therefore, $L_{m+i}V'=0$ for all $i>m$.
 Similarly, we can show that there exist
  $p,q\in\Z_{+}$ such that
 $J_{i}V'=G_{j}V'=0$ for all $i\geqslant p$, $j\geqslant q$.
 Hence, we get (5).

 $(5)\Rightarrow(6)$:
Fix $k,l,m\in\Z_{+}$ and a nonzero vector $v\in S$ such that $\tl^{(k,l,m)}v=0$.
By the PBW Theorem and the simplicity of $S$, $S$ has a spanning set
consisting of vectors of the form
$$G^{{\bf k}}J^{\bf i}L^{\bf j}v
=\cdots G_{k-2}^{k_{2}}G_{k-1}^{k_{1}}\cdots J_{l-2}^{i_{2}}J_{l-1}^{i_{1}}\cdots L_{m-2}^{j_{2}}L_{m-1}^{j_{1}}v,$$
where ${\bf k}=(\ldots,k_{2},k_{1}),{\bf i}=(\ldots,i_{2},i_1),{\bf j}=(\ldots,j_{2},j_{1})\in\M$.

{\bf Claim:}
$$L_{n}G^{{\bf k}}J^{\bf i}L^{\bf j}v=G_{n}G^{{\bf k}}J^{\bf i}L^{\bf j}v=J_{n}G^{{\bf k}}J^{\bf i}L^{\bf j}v=0$$
for $n> k+l+m+{\bf w(k+i+j)}$.

We prove this claim by induction on $s:={\bf d(k+i+j)}$.
The case of $s=0$ is clear.
For $s=1$, $G^{{\bf k}}J^{\bf i}L^{\bf j}v$ is of the form $G_{k-b}v$, $J_{l-b}v$, or $L_{m-b}v$
for some $b\in\N$.
Suppose $G^{{\bf k}}J^{\bf i}L^{\bf j}v=G_{k-b}v$.
For $n>k+l+m+b$, by assumption,
we have
\begin{eqnarray*}
\begin{aligned}
&L_{n}G_{k-b}v=G_{k-b}L_{n}v+[L_{n},G_{k-b}]v=0,\\
&J_{n}G_{k-b}v=G_{k-b}J_{n}v+[J_{n},G_{k-b}]v=0\;\;\mbox{and}\;\;G_{n}G_{k-b}v=0.
\end{aligned}
\end{eqnarray*}
Similarly, we can get the results for the cases of $ J_{l-b}v$ and $L_{m-b}v $.
Suppose $s>1$ and the claim holds for ${\bf d(k+i+j)}<s$.
Then for any $G^{{\bf k}}J^{\bf i}L^{\bf j}v$ with ${\bf d(k+i+j)}=s$,
$n> k+l+m+{\bf w(k+i+j)}$, we have
\begin{eqnarray}
 L_{n}G^{{\bf k}}J^{\bf i}L^{\bf j}v
&=&\sum\limits_{a,{\bf \tilde{k}}} G^{\bf \tilde{k}}[L_{n},G_{k-a}]G^{{\bf k-\tilde{k}}-\varepsilon_{a}}J^{\bf i}L^{\bf j}v
+\sum\limits_{a,{\bf \tilde{i}}} G^{\bf k}J^{\bf \tilde{i}}[L_{n},J_{l-a}]J^{{\bf i-\tilde{i}}-\varepsilon_{a}}L^{\bf j}v
\nonumber\\
&&{}
+\sum\limits_{a,{\bf \tilde{j}}} G^{\bf k}J^{\bf i}L^{\bf \tilde{j}}[L_{n},L_{m-a}]L^{{\bf j-\tilde{j}}-\varepsilon_{a}}v.
\end{eqnarray}
By induction hypothesis, we get $ L_{n}G^{{\bf k}}J^{\bf i}L^{\bf j}v=0$ for
$n> k+l+m+{\bf w(k+i+j)}$.
Similarly, we can show that $G_{n}G^{{\bf k}}J^{\bf i}L^{\bf j}v=J_{n}G^{{\bf k}}J^{\bf i}L^{\bf j}v=0$
for $n> k+l+m+{\bf w(k+i+j)}$.
Hence, $S$ is a restricted module.
This completes the proof.
\end{proof}

\bt{clasf1}
Suppose $F\ne 0$.
Let $S$ be a simple restricted $\tl$-module.
Assume that there exist $a,b\in\Z_{+}$ such that
the actions of $G_{a},L_{b}$ on $S$ are injective.
Then there exist $k,d,l\in\Z_{+}$ such that $k+d>0$ and
\vspace{3mm}\\
{\rm (i)} $N_{k,k+d,l}=\{ v\in S\ | \ G_{k+i}v=J_{k+d+i}v=L_{l+i}v=0 \,\,\mbox{for all}\,\,i\in\N\}\neq 0$; \\
{\rm (ii)} $G_{k},L_{l}$ act injectively on $ N_{k,k+d,l}.$
\vspace{3mm}\\
Furthermore, if $l>k+d$,
denote $V=N_{k,k+d,l}$,
then $V$ is a simple $\tl_{d}$-module
satisfying the conditions in Theorem \ref{idmod1}
and $S\cong {\rm Ind}_{\tl_{d}}^{\tl}(V)$.
\et
\begin{proof}
For any $r,s,t\in\Z$,
 consider the vector space
\begin{equation}\label{vsN}
N_{r,s,t}=\{ v\in S\ | \ G_{r+i}v=J_{s+i}v=L_{t+i}v=0 \,\,\mbox{for all}\,\,i\in\N\}.
\end{equation}
By Theorem \ref{eqchar} (5),
we know that $N_{r,s,t}\ne 0$ for sufficiently large integers $r,s,t$.
Note that if $N_{r,s,t}\ne 0$, then for any $i_{1},i_{2},i_{3}\in\Z_{+}$,
$N_{r+i_{1},s+i_{2},t+i_{3}}\neq 0$.
Since $G_{a},L_{b}$ act injectively on $S$, we can find
smallest integers $k\in\Z_{\geqslant a}$,
$l\in\Z_{\geqslant b}$ such
that $N_{k,k+d,l}\neq0$ for some $d\in\Z_{+}$ with $k+d>0$.
Assume that $G_{k}$ does not act injectively on $N_{k,k+d,l}$.
Then $N_{k-1,k+d,l}\neq 0$ which contradicts our choice of $k$.
Thus $G_{k}$ acts injectively on $N_{k,k+d,l}$.
Similarly, we have $L_{l}$ acts injectively on $N_{k,k+d,l}$.
Denote $V=N_{k,k+d,l}$. Then
$V\neq 0$.
If $l>k+d$,
it is straightforward to check that
$V$ is an $\tl_{d}$-module
satisfying the conditions in Theorem \ref{idmod1}.

There exists a canonical $\tl$-module epimorphism
\begin{equation}
\pi :\id_{\tl_d}^{\tl}(V) \rightarrow S,\;\;\pi (1\otimes v) =v\;\;\mbox{for any}\; v\in V.
\end{equation}
Let $K=\mbox{ker}(\pi)$ be the kernel of $\pi$.
 It is clear that $K\cap V=0$.
 If $K\ne 0$, we can choose a nonzero vector $v\in K \backslash V$
  such that $\d(v)=({\bf k,i,j})$ is minimal possible with respect
  to the principal total order ${\bf o}$.
   Note that $K$ is an $\tl$-submodule of $\id_{\tl_d}^{\tl}(V)$
    and hence is stable under the actions of $L_i,J_i$ and $G_{i}$ for all $i\in\Z$.
   By Lemmas \ref{lemred1.1}, \ref{lemred1.2} and \ref{lemred1.3},
 we can create a new vector $u\in K$ with $\d(u)\prec ({\bf k,i,j})$,
 which is a contradiction.
 Thus we have $K=0$, that is, $S\cong \id_{\tl_d}^{\tl}(V) $.
 By the property of induced modules, we know $V$ is a simple $\tl_d$-module.
\end{proof}

\bt{clasf2}
If $F=0$ and $C_{3}\neq 0$.
Let $S$ be a simple restricted $\tl$-module.
Assume that there exists $a\in \Z_{+}$ such that the action of $G_a$ on $S$ is injective.
Then there exists a simple $\tl_{\ud}$-module $V$ for $\ud\in\Z_{+}^{2}$
satisfying the conditions in Theorem \ref{idmod2} such that $S\cong {\rm Ind}_{\tl_{\ud}}^{\tl}(V)$.
\et
\begin{proof}
 For $r,s,t\in\Z$, consider the vector space $N_{r,s,t}$ in (\ref{vsN}).
Then $N_{r,s,t}\ne 0$ for sufficiently large integers $r,s,t$.
Since $G_{a}$ acts injectively on $S$, we can find a
smallest integer $k\in\Z_{\geqslant a}$ such
that $N_{k,s,t}\neq0$ for some $t\geqslant 2s\geqslant k$.
Assume that $G_{k}$ does not act injectively on $N_{k,s,t}$.
Then $N_{k-1,s,t}\neq 0$ which contradicts to our choice of $k$.
Thus $G_{k}$ acts injectively on $N_{k,s,t}$.
Denote $d_{1}=t-k$, $d_{2}=s$, $\ud=(d_{1},d_{2})$ and $V=N_{k,d_{2},k+d_{1}}\neq 0$.
Then
$V$ is an $\tl_{\ud}$-module
satisfying the conditions in Theorem \ref{idmod2}.

Since $S$ is a simple $\tl$-module,
 there exists a canonical $\tl$-module epimorphism
\begin{equation}
\pi :\id_{\tl_{\ud}}^{\tl}(V) \rightarrow S,\;\;\pi (1\otimes v) =v,\;\;\mbox{for any}\; v\in V.
\end{equation}
It follows from Lemma \ref{lemred2}
that $\mbox{Ker}\;\pi=0$.
Hence $S\cong \id_{\tl_{\ud}}^{\tl}(V)$,
and $V$ is a simple $\tl_{\ud}$-module
due to the property of induced modules.
\end{proof}

\justifying

\end{document}